\documentclass[12pt,a4paper,twoside,final,notitlepage,leqno]{article}
\usepackage[english]{babel}
\usepackage[T1]{fontenc}
\usepackage{epsfig, graphicx, amssymb}
\usepackage{amsmath,amsthm,epsfig,amsfonts}
\usepackage{float}
\usepackage{color}
\newcommand{\monitem}{ \smallskip \noindent $\bullet$ \quad  }
\newcommand{\moneq}{\vspace*{-6pt} \begin{equation} \displaystyle }
\newcommand{\moneqstar}{\vspace*{-6pt} \begin{equation*} \displaystyle }
\newcommand{\monendstar}{\vspace*{-6pt} \end{equation*}   }
\newcommand{\monend}{\vspace*{-6pt} \end{equation}   }
\newcommand{\moneqarraystar}{ \begin{eqnarray*} \displaystyle }
\newcommand{\monendarraystar}{ \end{eqnarray*}   }

\newcommand{\dd}{{\rm d}}
\newcommand{\R}{\mathbb{R}}

\definecolor{vertfonce}{rgb}{0.0, 0.5, 0.0}

\def\section*#1{}

\usepackage{fancyhdr}
\renewcommand{\headrulewidth}{0pt}

\begin{document}

\fancypagestyle{plain}{ \fancyfoot{} \renewcommand{\footrulewidth}{0pt}}
\fancypagestyle{plain}{ \fancyhead{} \renewcommand{\headrulewidth}{0pt}}

~

  \vskip 2.1 cm

\centerline {\bf \LARGE A  lattice Boltzmann scheme with equilateral}

 \bigskip 

\centerline {\bf \LARGE triangles for diffusion and acoustics}

\bigskip  \bigskip \bigskip

\centerline { \large   Fran\c{c}ois Dubois$^{ab}$ and Pierre Lallemand$^{c}$}

\smallskip  \bigskip

\centerline { \it  \small
  $^a$ Laboratoire de Math\'ematiques d'Orsay, Facult\'e des Sciences d'Orsay,}

\centerline { \it  \small   Universit\'e Paris-Saclay, France.}

\centerline { \it  \small
$^b$    Conservatoire National des Arts et M\'etiers, LMSSC laboratory,  Paris, France.}

\centerline { \it  \small
  $^c$   Beijing Computational Science Research Center}

\centerline { \it  \small 
  East Xibeiwang Road, Haidian District, Beijing, China.}

\bigskip  \bigskip  \bigskip

\centerline {{25 June 2026}
{\footnote {\rm  \small $\,$ This contribution will be submitted for publication shortly.}}}




\bigskip  \bigskip
\noindent {\bf \large Abstract}

\noindent
This contribution studies the Boltzmann scheme on a ``D2T4''grid constructed on meshes using equilateral triangles.
The center of each triangle is connected to itself and to three other triangles {\it via} the edges of the mesh.
We adopt the multiple relaxation time approach. 
Applications for diffusion and acoustics problems are considered.
Consistency analysis is particularly delicate. We propose an approach based on taking bipoints into account.
We derive  equivalent partial differential equations for diffusion and acoustics. 
These  systems of equations  are then approximated numerically
using the D2T4 lattice Boltzmann method.
A comparison with an analytical calculation in the case of periodic boundary conditions shows the convergence of the
D2T4 lattice Boltzmann scheme. 

\noindent

\newpage

\bigskip \bigskip    \noindent {\bf \large    1) \quad  Introduction}

\fancyhead[EC]{\sc{Fran\c{c}ois Dubois and Pierre Lallemand}}
\fancyhead[OC]{\sc{A  lattice Boltzmann scheme with equilateral triangles}}

\fancyfoot[C]{\oldstylenums{\thepage}}

\smallskip \noindent
The very first gas network of Hardy, Pomeau and Pazzis \cite{HPP73} used a square grid. 
Due to some defects in the asymptotic analysis, 
Frisch, Hasslacher, and  Pomeau \cite {FHP86}
propose to formulate cellular automata on equilateral triangles.
With modern notations, the corresponding lattice Boltzmann scheme
could be named as ``D2T7''.
After this success an  extension to three spatial dimensions was proposed by 
 d'Humi\`eres {\it et al.} \cite{HLF86}. 
However, the development of the mesoscopic approach of lattice Boltzmann schemes
with  
Higuera and Jim\'enez \cite{HJ89}, 
Higuera, Succi and Benzi \cite{HSB89},
Qian, d'Humi\`eres and Lallemand \cite {QHL92} and many others, 
enabled the formulation of the model on square and cubic grids.

\begin{figure}    [H]  \centering
\centerline{\includegraphics[width=.45\textwidth]{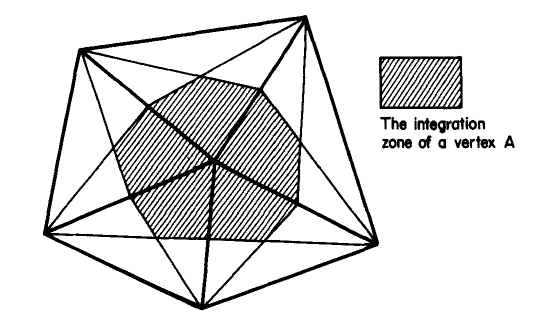}} 
\vskip -.5 cm 
\caption{Control volume around a vertex of a triangular mesh ; figure extracted from Angrand and Dervieux \cite{AD84}.}
\label{volume-fini-Dervieux} \end{figure}

\smallskip \noindent
Completely independently, the issue of calculating a gradient at a vertex 
for triangular grids with linear finite elements  was addressed
in the context of finite volumes by Angrand, Dervieux  {\it et al.} \cite {ADBPV83,AD84}. 
These authors propose a construction of a finite volume around  vertices of a general conforming  triangular mesh,
as illustrated in Figure \ref{volume-fini-Dervieux}.

\vskip -.3 cm 

\noindent 
This type of discrete geometry has been used by  Chen \cite {Ch98} 
for his volumetric formulation of lattice Boltzmann schemes on irregular grids.
We refer also to Peng {\it et al.} 
 \cite{PXDC99} for such an approach. 
Karlin, Succi and Orszag, \cite {KSO99}
consider a class of lattices which are structured in the sense that the number of links per lattice site (the
connectivity number) is constant throughout the lattice.
The cell-vertex finite volume discretization extended to three space dimensions in
Rossi {\it et al.}  
\cite{RUBS05}. 
Moreover,   Pontrelli {\it et al.} 
\cite{PUS09}   
have extended this formulation for non-Newtonian ﬂows.
The volumetric formulation has been generalized 
  to arbitrary coordinates by Chen \cite{Ch21}. 
However, this approach requires an interpolation phase that introduces numerical viscosity
into the lattice Boltzmann approach, similar to the Lax-Wendroff scheme \cite{LW60}.
One of the fundamental aspects of the Boltzmann lattice approach is lost in favor
of greater flexibility in the meshing process.

\smallskip \noindent
In his work on the diffusion equation for triangular meshes, 
van der Sman \cite{VdS04}
is using Voronoi cells constructed from a family of Delaunay triangulation.  
His approach for very general meshes does not introduce interpolation and  
we believe his results are fundamental.
In our contribution  \cite {DL13},  
we resume our study of lattice Boltzmann scheme on equilateral triangles.
First with the D2T7 model, which places degrees of freedom at the vertices of the mesh,
and with the D2T4  model, which places the calculation nodes at the center of gravity
of the triangles. 
We had limited ourselves to isotropic diffusion processes,
and the analysis we had done of the scheme had some flaws. 
In this work, we revisited this study on the D2T4 stencil
using equilateral triangles. On the one hand, with a view to applications
in thermal engineering, but also for acoustics. 


\smallskip \noindent
The  outline of our contribution is the following.
In Section~2, we describe the D2T4 scheme on equilareral triangles,
ensuring that, when necessary, the diffusion case is separated from the acoustic case.
For the asymptotic analysis, we remark that 
the dynamics of the lattice  operates on a two-point stencil. 
In order to treat it in as conventional a manner as possible,
we introduce the concept of breathing moments in Section~3. 
Then in Section 4, we derive the equivalent partial differential equations
for diffusion  and acoustic models.
Numerical results for periodic flows are presented in Section~5.
Some words of conclusion are proposed in Section 6. 
Three appendices present formal calculations based on a
decoupling of the two families of triangles.

\bigskip \bigskip    \noindent {\bf \large    2) \quad  {Description of the D2T4 lattice Boltzmann scheme}}

\smallskip \noindent
In this section, we recall the basics about the  D2T4 lattice Boltzmann scheme
for equilateral triangles.
We are essentially following the outline described in our contribution \cite{DL13},
although the notations used here are slightly different.

\vskip -1.3 cm 
\begin{figure}    [H]  \centering
\centerline{\includegraphics[width=.77 \textwidth]{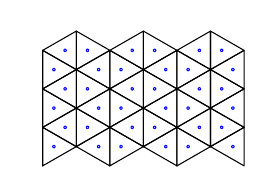}}
\vskip -1.3 cm 
\caption{Mesh with equilateral triangles; the degrees of freedom of the D2T4 algorithm
are the centers of the triangles, represented by small blue circles.}
\label{maillage-triangles} \end{figure}
\bigskip 
\vskip -.5 cm

\monitem {\bf  Discrete velocities}

We consider a regular mesh composed with equilateral triangles,
as the one of Figure~\ref{maillage-triangles} to fix the ideas. 
Unlike the initial choice made by Frisch, Haslacher, and Pomeau \cite{FHP86}, who place the particles
at the vertices of the mesh,
we assume that the particles are placed at the centers of the triangles.
With this framework, each triangle has three neighbors through its three edges.
We suppose also that is is connected to itself
with the null velocity. We denote  $ \, \Delta x \, $ the distance between the centers of two neiboring triangles. 
Then four velocities connect a given triangle with its neighbors, which justifies the name ``D2T4''.
In the example of  Figure~\ref{maillage-triangles}, any  vertical edge  separates two triangles: one on the left, named
$\, x^\ell \,$ and one on the right, $ \, x^r $.
With this convention, each triangle of the mesh belongs to one of the two familes, 
on left-type and right-type triangles. 

\begin{figure}    [H]  \centering
\centerline{\includegraphics[width=.55 \textwidth]{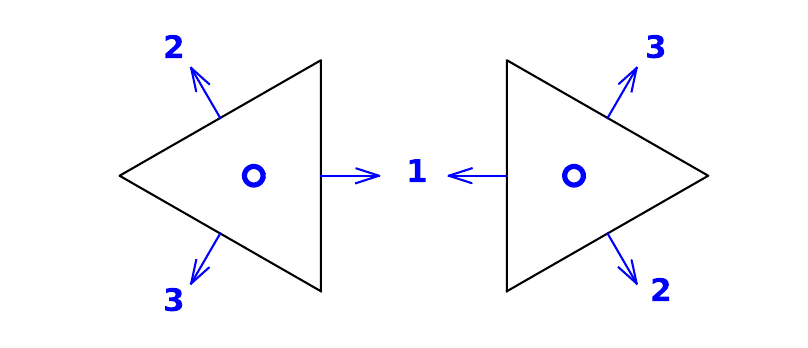}} 
\vskip -.2 cm 
\caption{Two facing triangles $ \, x^\ell \, $ of left-type and $ \, x^r \, $ of right type.}
\label{fig-bi-triangle} \end{figure}

\smallskip \noindent
We introduce  a time state $ \, \Delta t $. 
A reference scale speed~$ \, \lambda \,$ is associated to the ratio between the space step $ \, \Delta x \, $
and the time step $ \, \Delta t $:
\moneqstar 
 \lambda = {{\Delta x}\over{\Delta t}}  \, .
 \monendstar
\smallskip \noindent
We consider two families $ \, v_j^\ell \, $ and $ \, v_j^r \, $ of discrete velocities, for $ \, 0 \leq j \leq 3 \, $
for  the D2T4 scheme described in Figure~\ref{fig-bi-triangle}.
The first family of velocities  $ \, v_j^\ell \, $ go from the triangle $ \, x^\ell \, $ towards the four neighboring triangles; we
have
\moneq \label{vitesses-gauche}
\{   v_j^\ell \} = \begin{pmatrix}
0 & \lambda &  -{1\over2}\, \lambda & {1\over2}\, \lambda  \\
0 & 0 &  {{\sqrt{3}}\over2}\, \lambda & -{{\sqrt{3}}\over2}\, \lambda  \end{pmatrix} \, . 
\monend
The numbering corresponds to the convention presented in  Figure~\ref{fig-bi-triangle}.
 \vskip -.5 cm

\smallskip 
We note that the zero neighbor of a left-type triangle is of the same type, while the other three neighbors are right-type. 
The property is analogous for right-type triangles. The neighbor with number zero is of the same type,
while the other three neighbors are left-type triangles. 
The second family family of velocities  $ \, v_j^r \, $ go from the triangle $ \, x^r  \, $ towards its neighbors.
We have the simple relation 
\moneq \label{vitesses-droite}
 v_j^r +  v_j^\ell = 0 \,, \quad 0 \leq j \leq 3 \, 
\monend 
with a numbering convention proposed in  Figure~\ref{fig-bi-triangle}.
The four neighbors $ \, x^\ell_j \, $   of a left-type triangle  $ \, x^\ell \, $ can be written as
\moneqstar
x^\ell_j  =  x^\ell + \Delta t \, v_j^\ell \,, \quad 0 \leq j \leq 3 \, .
\monendstar 
The neighbor  $ \, x^\ell_0 \, $ number zero of the triangle  $ \, x^\ell \, $ is  simply  $ \, x^\ell \, $ itself.
It is a left-type triangle.
The three other neighbors belong to the family of the right-type triangles.
It is clear when looking to  Figure~\ref{maillage-triangles}.
Analogously, the four neighbors $ \, x^r_j \, $   of a right-type triangle  $ \, x^r \, $ admit the expression 
\moneqstar
x^r_j  =  x^r + \Delta t \, v_j^r \,, \quad 0 \leq j \leq 3 \, .
\monendstar 
The triangle $ \, x^r \, $ is  its own neighbor  $ \, x^r_0 \, $ number zero, and belongs to the right family.
The three other neighbors  $ \, x^r_j \, $  for $\, 1 \leq j \leq 3 \, $ belong to the family of  left triangles.

\bigskip \monitem {\bf  Particles and moments}

\noindent
A cell center of the triangular lattice lattice $\, {\cal L} \, $ is denoted by $\, x $.
In other words, the position~$ \, x \, $ located at the barycentre of a triangle of the mesh.
The ``outgoing'' particles $ \, f(x) \, $ join a given node $ \, x \, $ towards its
four neighbors. Conversely, the ``incoming'' particles  $ \, g(x) \, $ join
the neighbors of the  triangle $ \, x \, $ to itself.

\smallskip \noindent 
We must be mindful that triangles labeled ``left'' may have different particle distributions
than triangles labeled ``right''. 
We therefore adopt the notation  $\, f_\ell(x^\ell) \, $ (respectively $\, f_r(x^r) $) 
for particles leaving the  triangle $ \, x^\ell \,$ (respectively leaving the triangle $ \, x^r $)
and  $\, g_\ell(x^\ell) \, $ (respectively $\, g_r(x^r) $) for particles entering the  triangle
$ \, x^\ell \,$ (respectively entering the triangle~$ \, x^r $). 

\smallskip \noindent 
We define the local moments in each triangle to connect the two notions. We introduce
a ``left'' matrix of moments $\, M^\ell \,$ with the discrete distribution of particles,
in the spirit proposed by d'Humi\`eres \cite{DDH92}:
\moneq  \label{matrice-M-gauche}
M^\ell  =  \begin{pmatrix} 1 & 1 & 1 & 1 \\
0 & \lambda &  -{1\over2}\, \lambda & {1\over2}\, \lambda  \\
0 & 0 &  {{\sqrt{3}}\over2}\, \lambda & -{{\sqrt{3}}\over2}\, \lambda   \\
-3 \, \lambda^2 &  \lambda^2 &  \lambda^2 &  \lambda^2 \end{pmatrix} \, . 
\monend
The first line define the moment $ \, \rho \, $ named ``density''. The second and third lines of the
matrix~$ \, M^\ell \, $ defined in (\ref{matrice-M-gauche}) define
the two components $ \, J_x \, $ and $ \, J_y \, $ of the momentum.
They reproduce the components of the four velocities $ \, v^\ell_j $.
The  last line introduces the energy $\, e $. We observe that the matrix  $ \, M^\ell \, $
introduced in (\ref{matrice-M-gauche}) has orthogonal lines:
\moneqstar
\sum_{j=0}^3  \big(M^\ell\big)_{ij} \, \big(M^\ell\big)_{kj}   =  0 \, \quad {\rm for} \,\,\,  i \not= k \, . 
\monendstar 
The moments $ \, m(x^\ell) \equiv \big( \rho^\ell ,\, J_x^\ell ,\, J_y^\ell  ,\, e^\ell  \big)^{\rm t} \, $ of a left triangle
are defined according to a classical relation
\moneq  \label{f-to-m-gauche}
 m(x^\ell)  = M^\ell  \,  f_\ell(x^\ell) \,. 
\monend
For a right triangle, we just change the signs of the velocities:
\moneq  \label{matrice-M-droite}
M^r  =  \begin{pmatrix} 1 & 1 & 1 & 1 \\
0 & -\lambda &  {1\over2}\, \lambda & -{1\over2}\, \lambda  \\
0 & 0 &  -{{\sqrt{3}}\over2}\, \lambda & {{\sqrt{3}}\over2}\, \lambda   \\
-3 \, \lambda^2 &  \lambda^2 &  \lambda^2 &  \lambda^2 \end{pmatrix} \, . 
\monend
We have again an orthogonality relation 
\moneqstar
\sum_{j=0}^3  \big(M^r\big)_{ij} \, \big(M^r \big)_{kj}   =  0 \, \quad {\rm for} \,  i \not= k \, . 
\monendstar 
The moments $ \, m(x^r) \equiv \big( \rho ,\, J_x ,\, J_y ,\, e \big)^{\rm t} \, $ of a rith triangle
are defined according to
\moneq  \label{f-to-m-droite}
 m(x^r)  = M^r  \,  f_r(x^r) \,. 
\monend 
With the incoming particles $\, g $, the sign of the velocities is inverted. We have  
\moneq  \label{g-to-m} \left\{    \begin{array} {l}
 m(x^\ell)  = M^r  \,  g_\ell(x^\ell) \\
 m(x^r)  = M^\ell  \,  g_r(x^r) \,. \end{array} \right.  
\monend 

\smallskip \noindent
Although the distribution of particles is {\it a priori} discontinuous when moving from
a left-type triangle to a right-type triangle, 
moments ultimately represent real physical quantities and are supposed to be continuous functions of space and time.
So we do not use the labels ``$\ell$`` or ``$ r $''  for moments, simply denoted by $ m $,
as in relations  (\ref{f-to-m-gauche}), (\ref{f-to-m-droite}), and (\ref{g-to-m}).

\smallskip \noindent
We call ``breathing of the lattice'' the change of representation between ingoing and outgoing particles.
With the notations $ \, M^{-\ell} \equiv (M^\ell)^{-1} \, $ and $ \, M^{-r} \equiv (M^r)^{-1} $, we have
%
%
\moneqstar \left\{    \begin{array} {l}
 g_\ell(x^\ell) =   M^{-r} \,  M^{\ell}  \, f_\ell (x^\ell) =  M^{-\ell} \,  M^{r}  \, f_r (x^\ell)  \\ 
 g_r(x^r) =   M^{-\ell} \,  M^{r}   \, f_r(x^r) =   M^{-r} \,  M^{\ell}  \, f_r(x^r)  \,,  \end{array} \right.  \monendstar 
with
\moneqstar
M^{-r} \,  M^{\ell}  =  M^{-\ell} \,  M^{r}  =  \begin{pmatrix} 1 & 0 & 0 & 0 \\   0 & -{1\over3} &  {2\over3} & {2\over3} \\
  0  &  {2\over3} & -{1\over3} & {2\over3} \\  0  &  {2\over3} &  {2\over3} & -{1\over3} \end{pmatrix} \, . 
\monendstar

\smallskip \noindent 
We observe that the breathing of the lattice takes the same algebraic expression for left and right-type triangles. 

\bigskip \monitem {\bf  Relaxation of the moments}

\noindent
We consider in this contribution diffusion and acoustic  problems.
In the first case, only the first moment (density) is conserved whereas
density and the two components of momentum are conserved for acoustics. We set
in the diffusion case 
\moneq \label{WY-thermique}  \left\{    \begin{array} {rlcrl}
W^\ell =& \!\!\! (\rho^\ell) \,, &&  Y^\ell =& \!\!\! \big( J_x^\ell ,\, J_y^\ell  ,\, e^\ell  \big)^{\rm t} \\ 
W^r =& \!\!\! (\rho^r) \,, &&  Y^r =& \!\!\! \big( J_x^r ,\, J_y^r  ,\, e^r  \big)^{\rm t} \, . 
\end{array} \right. \monend 
This decomposition explicit the two families of moments; the conserved moments $ \, W \, $ and the non-cnserved moments $ \, Y $.
We have 
\moneq \label{moments-gd-WY}  
m(x^\ell) = \begin{pmatrix} W^\ell \\  Y^\ell \end{pmatrix} \,\, {\rm and} \quad 
m(x^r) = \begin{pmatrix} W^r \\  Y^r \end{pmatrix} \, . 
\monend 
The decomposition (\ref{moments-gd-WY}) is still valid in the acoustic case, except that the
definition of the conserved moments $\, W \, $ and the non-conserved ones $ \, Y \, $ has to be modified.
We set for acoustics 
\moneq \label{WY-acoustique}  \left\{    \begin{array} {rlcrl}
W^\ell =& \!\!\!  \big( \rho^\ell  ,\, J_x^\ell ,\, J_y^\ell \big)^{\rm t}   \,, &&  Y^\ell =& \!\!\! (e^\ell )   \\ 
W^r =& \!\!\!  \big( \rho^r  ,\, J_x^r ,\, J_y^r \big)^{\rm t}   \,, &&  Y^r =& \!\!\! (e^r )  \, . 
\end{array} \right. \monend 
The vector $ \, \Phi \, $ of equilibria is a function of the conserved moments:
\moneqstar
Y^{\rm eq} \equiv \Phi(W) \, .
\monendstar 
In the following, we consider only linear equilibria parameterized by a given scalar coefficient~$ \, \alpha $. 
In the pure diffusive  case, we have
\moneq \label{equilibre-thermique} 
\Phi(W) = \Phi(\rho) = \begin{pmatrix} 0 \\ 0 \\ \alpha \, \lambda^2 \, \rho \end{pmatrix} \, . 
\monend
Then
\moneqstar
Y^{\ell\, {\rm eq}} = \big( 0 ,\, 0 ,\, \alpha \, \lambda^2 \, \rho^\ell \big)^{\rm t} \,, \quad 
Y^{r, {\rm eq}} = \big( 0 ,\, 0 ,\, \alpha \, \lambda^2 \, \rho^r \big)^{\rm t} \, . 
\monendstar
The moments after relaxation follow a general relation introduced by d'Humi\`eres \cite{DDH92}. 
\moneqstar
m^* = \begin{pmatrix} W^* \\ Y^* \end{pmatrix} = \begin{pmatrix} W \\ Y + S \, (\Phi(W) - Y) \end{pmatrix} \, . 
\monendstar
The square matrix $ \, S \, $ is supposed to be a constant and diagonal in the conventional relation.
In the diffusion  case, we have
\moneqstar 
S = {\rm diag} \big( s_j ,\, s_j, \, s_e \big) \, .
\monendstar 
The relaxation of non-conserved moments follow the relations
\moneqstar \left\{    \begin{array} {lll}
  J_x^{\ell\, *} = (1-s_j) \, J_x^\ell  \,,\quad &  
  J_y^{\ell\, *} =  (1-s_j) \, J_y^\ell  \,,\quad &  
  e^{\ell\, *} =  (1-s_e) \, e^\ell + s_e \, \alpha \, \lambda^2 \, \rho^\ell \\ 
  J_x^{r\, *} = (1-s_j) \, J_x^r  \,,\quad &  
  J_y^{r\, *} =  (1-s_j) \, J_y^r  \,,\quad &  
  e^{r\, *} =  (1-s_e) \, e^r + s_e \, \alpha \, \lambda^2 \, \rho^r \,. 
\end{array} \right. \monendstar
Joined with the relations $  \,  \rho^{\ell\, *} =  \rho^{\ell} \, $ and $ \,  \rho^{r\, *} =  \rho^{r} $, the
relaxations of the moments for the left and right triangles takes the form
\moneq \label{relaxation-moments} 
 m^*(x^\ell) = J_0 \,\, m(x^\ell) \,,\quad m^*(x^r) = J_0 \,\, m(x^r)
\monend
with
\moneq \label{J0-thermique} 
J_0 = \begin{pmatrix} 1 & 0 & 0 & 0 \\ 0 & 1-s_j & 0 & 0 \\ 0 & 0 & 1-s_j & 0 \\
  \alpha \, \lambda^2 \, s_e & 0 & 0 & 1-s_e \end{pmatrix} \, . 
\monend 
For acoustics, the relations (\ref{equilibre-thermique}) to (\ref{J0-thermique})  take
an other algebraic form, to take into account the conserved moments (\ref{WY-acoustique}). We have now
\moneq \label{equilibre-acoustique} 
\Phi(W) = \Phi(\rho ,\, J_x, \, J_y) = \big( \alpha \, \lambda^2 \, \rho \big) \,,  
\monend

\vskip -.4 cm 
\moneqstar 
S = {\rm diag} \big( s_e \big) \, , 
\monendstar 
and the relaxation of the last non-conserved moment is simply 
\moneqstar \left\{    \begin{array} {l}
  e^{\ell\, *} =  (1-s_e) \, e^\ell + s_e \, \alpha \, \lambda^2 \, \rho^\ell \\ 
  e^{r\, *} =  (1-s_e) \, e^r + s_e \, \alpha \, \lambda^2 \, \rho^r \,. 
\end{array} \right. \monendstar
The relation (\ref{relaxation-moments}) is still valid, with a slight different matrix $ \, J_0 $:
\moneq \label{J0-acoustique} 
J_0 = \begin{pmatrix} 1 & 0 & 0 & 0 \\ 0 & 1 & 0 & 0 \\ 0 & 0 & 1 & 0 \\
  \alpha \, \lambda^2 \, s_e & 0 & 0 & 1-s_e \end{pmatrix} \, . 
\monend 
In a way, it is sufficient to set  $ \, s_j=0 \, $ to go from thermics to  acoustics. 

\bigskip \monitem {\bf  Free transport of the particles}

\noindent
Once the moments after relaxation $ \, m^*(x^\ell) \, $ and $ \, m^*(x^r) \, $ have been evaluated
in the triangles~$ \, x^\ell \, $ and $ \, x^r $,   
the two families of triangles exchange information during the advection step  
between the times $ \, t \, $ and $ \, t + \Delta t $.
We first go back to the particle distribution, inverting the relations
(\ref{f-to-m-gauche}) and (\ref{f-to-m-droite}): 
\moneqstar
f_\ell^*(x^\ell) = M^{-\ell} \, m^*(x^\ell) \,, \quad f_r^*(x^r) = M^{-r} \, m^*(x^r) \, . 
\monendstar 

\smallskip \noindent
The advection scheme must be specified for zero velocity on the one hand and for non-zero velocities on the other. 
First, each triangle~$ \, x^\ell \, $ exchange  information with itself with the help
of the particles of zero velocity. The particle after relaxation  $ \, f_{\ell 0}^* (x_\ell) \, $ is entering
into the same triangle at the new time step: 
\moneq \label{g-gauche-zero}
g_{\ell 0}(x^\ell,\, t+\Delta t) = f_{\ell 0}^* (x^\ell, \, t) \, . 
\monend
In a similar way,  each triangle~$ \, x^r \, $ is coupled with itself with the particle of zero velocity:
\moneq \label{g-droite-zero}
g_{r 0}(x^r,\, t+\Delta t) = f_{r 0}^* (x^r , \, t) \, . 
\monend
 \vskip -.5 cm
\begin{figure}    [H]  \centering
\centerline{\includegraphics[width=.77 \textwidth]{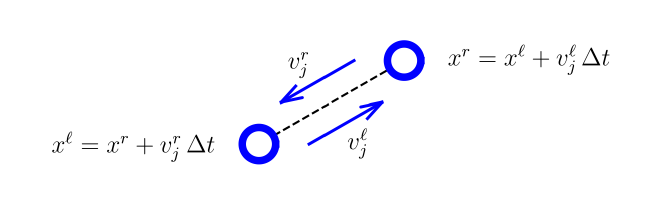}}
\vskip -.4 cm 
\caption{Exchange of particles between the two neighboring triangles  $ \, x^\ell \, $ and $ \, x^r $.}
\label{fig-bi-vitesses} \end{figure}

\smallskip 

%
\noindent 
Secondly, each triangle exchange with its three neighbors with the velocities numbered 1, 2, and 3.
The processus is decribed in Figure~\ref{fig-bi-vitesses}.
From left to right, the particle $ \, f_{\ell \,j}^* (x^\ell, \, t) \,$ go  outside the triangle $ \, x^\ell \, $
at time $ \, t \, $ and is entering at time $ \, t + \Delta t \, $ into  a right triangle with the neighboring number 
 $ \, j $: $ \, x_j^\ell =  x^r = x^\ell + v_j^\ell \, \Delta t $. We have 
\moneqstar
g_{r j}(x^\ell   + v_j^\ell \, \Delta t ,\, t+\Delta t) = f_{\ell j}^* ( x^\ell , \, t) \, .
\monendstar 
We exchange the role of the two triangles $ \, x^\ell \, $ and $ \, x^r $. Then
\moneqstar
g_{rj}(x^r ,\, t+\Delta t) = f_{lj}^* ( x^r - v_j^\ell \, \Delta t  , \, t)
=  f_{lj}^* ( x^r + v_j^d \, \Delta t  , \, t) =  f_{lj}^* ( x_j^r ,\, t ) \,. 
\monendstar 
We have established the relation 
\moneq \label{g-droite-j}
g_{rj}(x^r ,\, t+\Delta t) =  f_{lj}^* ( x_j^r ,\, t ) \, ,\quad 1 \leq j \leq 3 \, . 
\monend
From right to left, the particle $ \, f_{r j}^* (x^r , \, t) \,$ is going  outside the triangle $ \, x^r \, $
at time $ \, t \, $ and is entering at time $ \, t + \Delta t \, $ inside the left neighbor triangle with 
number $ \, j $: $ \, x_j^r =  x^\ell = x^r + v_j^r \, \Delta t $. We have 
\moneqstar
g_{lj}(x^r + v_j^r \, \Delta t ,\, t+\Delta t) = f_{rj}^* ( x^r , \, t) \, .
\monendstar 
We write this relation under the form
\moneqstar
 g_{lj}(x^\ell ,\, t+\Delta t) = f_{rj}^* ( x^\ell - v_j^r \, \Delta t  , \, t)
 =  f_{rj}^* ( x^\ell + v_j^\ell \, \Delta t  , \, t) =  f_{rj}^* ( x_j^\ell ,\, t) \, . 
\monendstar 
Then 
\moneq \label{g-gauche-j}
g_{\ell j}(x^\ell ,\, t+\Delta t) = f_{r j}^* (x_j^\ell ,\, t) \, ,\quad 1 \leq j \leq 3 \, . 
\monend
We recover the distribution of moments at the new time using the relations~(\ref{g-to-m}).
The D2T4 scheme is now entirely defined.

\newpage 
\bigskip \bigskip    \noindent {\bf \large    3) \quad  {Breathing moments}}

\smallskip \noindent
Setting up equivalent equations is tricky for this D2T4 scheme.
For all conventional schemes such as D2Q9 \cite{QHL92, LL00}, D2Q13 \cite{Qi93},
D3Q19 \cite{HGKLL02} or D3Q27 (see {\it e.g.} \cite{FHHLPR87}), the degrees of freedom are aligned.
Even for the scheme with D2T7 triangles, the particle emission points align.
For the D2T4 scheme, the situation is different. A quick look at Figure~\ref{maillage-triangles} shows that
a particle cannot move in a straight line over several time steps.
The structure described in Figure~\ref{fig-bi-vitesses} provides all available information.
For two neighboring triangles, two particles exchange during a time step, without direct interaction
with the other degrees of freedom present in the mesh.
We noted in our previous work that the equivalent equation (relation (6.2) of \cite{DL13}) does not
fully correspond to numerical observations. In particular, an anisotropic term
of the type~$ \, (\partial_x^2 - 3 \, \partial_y^2)\, \rho \, $ 
from the third-order terms of a diffusion model onwards gives pause for thought.

\smallskip \noindent
In Appendix A, we detail the overly simplistic calculation that fails to take lattice breathing into account
in the diffusive  case. We obtain partial differential equations that are consistent between left-type
and right-type triangles, but only up to the second order. At third order, the partial differential
equations are different, which is not acceptable {\it a priori}.

\smallskip \noindent
In Appendix B, we repeat a similar study, but in the case of the acoustic system.
In this case, inconsistencies appear right from the first-order equation!

\smallskip \noindent
In this section, we revisit this study. The main idea is to replace a single node $ \, x \in {\cal L} \, $
by a bipoint $ \, X =(x^\ell,\, x^r) \,$ of two neighboring triangles of the mesh.
In a first proposition, we reformulate the D2T4 scheme described in the previous section under the form
\moneq \label{evolution-du22}
m (X,t+\Delta t) = \exp(-\Delta t\, \Lambda) \, m^* (X,t         ) \,, 
\monend
with $ \, m(X) \, $ a set of moments, 
$ \, \Lambda \, $ a differential operator representing the advection in the space of moments,
and  $ \, m^*(X) \, $ the family of  moments after relaxation.
After this first transformation, we use our general methodology described in \cite{Du22}
to obtain equivalent partial differential equations. 

\smallskip \noindent
From the previous relations, we introduce a family of decoupled moments for the bipoint $\, X \equiv (x^\ell,\, x^r) $.
We introduce a vector with 8 components by sticking two adjacent triangles together: 
\moneq \label{moments-decouples}
m_d(X) = \begin {pmatrix} m(x^\ell) \\ m(x^r) \end{pmatrix} \,. 
\monend 
We have also a velocity vector $ \, w \in (\R^2)^8 \, $ by combining the speeds $\, v_j^\ell \, $ and $ \, v_j^r \, $
defined with the relations  (\ref{vitesses-gauche}) and (\ref{vitesses-droite}):
\moneq \label{vitesses-globales}
w = \big( v_0^\ell ,\, 0 ,\, v_2^\ell - v_1^\ell ,\, v_3^\ell - v_1^\ell  ,\, v_0^r ,\, 0 ,\, v_2^r -v_1^r  ,\, v_3^r -v_1^r \big) \,. 
\monend 
With this velocity dynamic, a left-type triangle is always related to another
triangle of the same type. This observation also applies to right-type triangles.

\bigskip 
{\bf  Proposition 1. Discrete evolution of a bipoint}

With the notations   (\ref{moments-decouples}) and (\ref{vitesses-globales}),  the D2T4 scheme defined at
section~2 can be written
\moneq \label{evolution-moments-decouples}
m_d(X,\, t+\Delta t)  = \begin{pmatrix} M^r & 0 \\ 0 & M^\ell \end{pmatrix} \,\, \Pi_\sigma \,\,
\exp (-\Delta t \,\, w.\nabla) \,\, \begin{pmatrix} M^{-\ell} & 0 \\ 0 & M^{-r}  \end{pmatrix} \, m_d^*(X,\, t) \,,
\monend
with
 $ \, m_d^*(X) \equiv \big( m^*(x^\ell) ,\, m^*(x^r) \big)^{\rm t} $, 
$ \, M^\ell \, $ and $ \, M^r \, $ defined in (\ref{matrice-M-gauche}) and (\ref{matrice-M-droite}),
a discrete permutation operator $ \, \Pi_\sigma \, $ defined by 
\moneq \label{Pi-sigma} 
\Pi_\sigma = \begin{pmatrix} 1 & 0 & 0 & 0 & 0 & 0 & 0 & 0 \\  0 & 0 & 0 & 0 & 0 & 1 & 0 & 0 \\
  0 & 0 & 0 & 0 & 0 & 0 & 1 & 0 \\  0 & 0 & 0 & 0 & 0 & 0 & 0 & 1   \\ 0 & 0 & 0 & 0 & 1 & 0 & 0 & 0  \\
 0 & 1 & 0 & 0 & 0 & 0 & 0 & 0 \\ 0 & 0 & 1 & 0 & 0 & 0 & 0 & 0 \\ 0 & 0 & 0 & 1 & 0 & 0 & 0 & 0 
  \end{pmatrix} \,, 
\monend
the inverse matrices $ \, M^{-\ell} \equiv (M^\ell)^{-1} $, $ \, M^{-r} \equiv (M^r)^{-1} $, 
and  a diagonal~8 by 8 matrix $\, w.\nabla \, $
defined by the relations $ \,\, \xi.\nabla \equiv \xi_x \, \partial_x + \xi_y \, \partial_y \,\,$ and 
because $ \, v_0^\ell = v_0^r = 0 $, by 
\moneq \label{w.nabla} 
w.\nabla = {\rm diag} \, \big( 0,\, 0 ,\, (v_2^\ell - v_1^\ell) .\nabla ,\,  (v_3^\ell - v_1^\ell) .\nabla ,\,
0 ,\, 0 ,\, (v_2^r - v_1^r) .\nabla ,\, (v_3^r - v_1^r) .\nabla \big) \,.
\monend

\bigskip \monitem
 Proof of Proposition 1 

\smallskip \noindent
First, we must order the different particles associated to the relations (\ref{g-gauche-zero}) to (\ref{g-gauche-j}).
We have
\smallskip \noindent 
$ g_{\ell 0}(x^\ell ,\, t+\Delta t) = f_{\ell 0}^* (x^\ell ,\, t ) $   \qquad due to (\ref{g-gauche-zero}), 

\noindent 
$ g_{\ell 1}(x^\ell ,\, t+\Delta t) = f_{r 1}^* (x_1^\ell ,\, t) $ \qquad due to (\ref{g-gauche-j})

\noindent \qquad  \qquad  \qquad
$\,\,\,  = f_{r 1}^* (x^r  ,\, t) \, $ \qquad   because $ \, x_1^\ell = x^r $, 

\noindent 
$ g_{\ell 2}(x^\ell ,\, t+\Delta t) = f_{r 2}^* (x_2^\ell ,\, t) $ \qquad due to (\ref{g-gauche-j})

\noindent \qquad  \qquad  \qquad
$\,\,\,  = f_{r 2}^* (x^\ell + \Delta t \, v_2^\ell  ,\, t)  =  f_{r 2}^* (x^r + \Delta t \, v_1^r + \Delta t \, v_2^\ell  ,\, t)
= f_{r 2}^* (x^r - \Delta t \, ( v_2^r -  v_1^r)  ,\, t) $, 

$ g_{\ell 3}(x^\ell ,\, t+\Delta t) = f_{r 3}^* (x^r - \Delta t \, ( v_3^r -  v_1^r)  ,\, t) $ \qquad due to (\ref{g-gauche-j})
and the previous calculus,

\noindent 
$ g_{r 0}(x^r ,\, t+\Delta t) = f_{r 0}^* (x^r ,\, t ) $   \qquad due to (\ref{g-droite-zero}), 

\noindent 
$ g_{r 1}(x^r ,\, t+\Delta t) = f_{\ell 1}^* (x_1^r ,\, t) $ \qquad due to (\ref{g-droite-j})

\noindent \qquad  \qquad  \qquad
$\,\,\,  = f_{\ell 1}^* (x^\ell  ,\, t) \, $ \qquad   because $ \, x_1^r = x^\ell  $, 

\noindent 
$ g_{r 2}(x^r ,\, t+\Delta t) = f_{\ell 2}^* (x_2^r ,\, t) $ \qquad due to (\ref{g-droite-j})

\noindent \qquad  \qquad  \qquad
$\,\,\,  = f_{\ell 2}^* (x^r + \Delta t \, v_2^r  ,\, t)  =  f_{\ell 2}^* (x^\ell + \Delta t \, v_1^\ell + \Delta t \, v_2^r  ,\, t)
= f_{\ell 2}^* (x^\ell - \Delta t \, ( v_2^\ell -  v_1^\ell)  ,\, t) $, 

$ g_{r 3}(x^r ,\, t+\Delta t) = f_{\ell 3}^* (x^\ell - \Delta t \, ( v_3^\ell -  v_1^\ell)  ,\, t) $ \qquad due to (\ref{g-droite-j})
and the previous calculus.

\smallskip \noindent
Then we have
\moneqstar 
 \begin{pmatrix}  g_{\ell 0}(x^\ell ,\, t+\Delta t) \\ g_{\ell 1}(x^\ell ,\, t+\Delta t) \\ g_{\ell 2}(x^\ell ,\, t+\Delta t) \\ g_{\ell 3}(x^\ell ,\, t+\Delta t) \\ 
   g_{r 0}(x^r ,\, t+\Delta t) \\  g_{r 1}(x^r ,\, t+\Delta t) \\  g_{r 2}(x^r ,\, t+\Delta t) \\  g_{r 3}(x^r ,\, t+\Delta t)  \end{pmatrix}
 = \begin{pmatrix}  f_{\ell 0}^* (x^\ell ,\, t ) \\  f_{r 1}^* (x^r  ,\, t) \\  f_{r 2}^* (x^r - \Delta t \, ( v_2^r -  v_1^r)  ,\, t) \\ 
f_{r 3}^* (x^r - \Delta t \, ( v_3^r -  v_1^r)  ,\, t) \\  f_{r 0}^* (x^r ,\, t )  \\ f_{\ell 1}^* (x^\ell  ,\, t)  \\
f_{\ell 2}^* (x^\ell - \Delta t \, ( v_2^\ell -  v_1^\ell)  ,\, t) \\  f_{\ell 3}^* (x^\ell - \Delta t \, ( v_3^\ell -  v_1^\ell)  ,\, t) \end{pmatrix}
 = \Pi_\sigma \,\, 
\begin{pmatrix}  f_{\ell 0}^* (x^\ell ,\, t ) \\ f_{\ell 1}^* (x^\ell  ,\, t)  \\
f_{\ell 2}^* (x^\ell - \Delta t \, ( v_2^\ell -  v_1^\ell)  ,\, t) \\  f_{\ell 3}^* (x^\ell - \Delta t \, ( v_3^\ell -  v_1^\ell)  ,\, t) \\
 f_{r 0}^* (x^r ,\, t ) \\    f_{r 1}^* (x^r  ,\, t) \\  f_{r 2}^* (x^r - \Delta t \, ( v_2^r -  v_1^r)  ,\, t) \\ 
f_{r 3}^* (x^r - \Delta t \, ( v_3^r -  v_1^r)  ,\, t)  \end{pmatrix}
\monendstar
with the permutation matrix $ \, \Pi_\sigma \, $ defined in  (\ref{Pi-sigma}).
We observe now that we have
\moneqstar  \left\{    \begin{array} {rcl}
f_{\ell 0}^* (x^\ell ,\, t ) &=& \exp (-\Delta t \, w_0 . \nabla) \,  f_{\ell 0}^* (x^\ell ,\, t ) \\
f_{\ell 1}^* (x^\ell  ,\, t) &=& \exp (-\Delta t \, w_1 . \nabla) \,  f_{\ell 1}^* (x^\ell ,\, t ) \\
f_{\ell 2}^* (x^\ell - \Delta t \, ( v_2^\ell -  v_1^\ell)  ,\, t) &=& \exp (-\Delta t \, w_2 . \nabla) \,  f_{\ell 2}^* (x^\ell ,\, t ) \\
f_{\ell 3}^* (x^\ell - \Delta t \, ( v_3^\ell -  v_1^\ell)  ,\, t)  &=& \exp (-\Delta t \, w_3 . \nabla) \,  f_{\ell 3}^* (x^\ell ,\, t ) \\
f_{r 0}^* (x^r ,\, t ) &=& \exp (-\Delta t \, w_4 . \nabla) \,  f_{r 0}^* (x^r ,\, t ) \\
f_{r 1}^* (x^r  ,\, t) &=&  \exp (-\Delta t \, w_5 . \nabla) \,  f_{r 1}^* (x^r ,\, t ) \\
f_{r 2}^* (x^r - \Delta t \, ( v_2^r -  v_1^r)  ,\, t) &=& \exp (-\Delta t \, w_6 . \nabla) \,  f_{r 2}^* (x^r ,\, t ) \\
f_{r 3}^* (x^r - \Delta t \, ( v_3^r -  v_1^r)  ,\, t) &=& \exp (-\Delta t \, w_7 . \nabla) \,  f_{r 3}^* (x^r ,\, t ) \,.
\end{array} \right. \monendstar
We can write the previous relations in a compact form
\moneqstar 
\begin{pmatrix}  g_{\ell}(x^\ell ,\, t+\Delta t) \\   g_{r}(x^r ,\, t+\Delta t) \end{pmatrix} = 
\Pi_\sigma \,\, \exp (-\Delta t \, w . \nabla) \,\, \begin{pmatrix}  f_{\ell}^* (x^\ell ,\, t ) \\  f_{r}^* (x^r ,\, t ) \end{pmatrix} 
\monendstar
then due to the relations    (\ref{f-to-m-gauche}), (\ref{f-to-m-droite}) and      (\ref{g-to-m}), 
\moneqstar
m_d(X ,\, t+\Delta t) = \begin{pmatrix} M^r & 0 \\ 0 & M^\ell \end{pmatrix} \, \Pi_\sigma \,\, \exp (-\Delta t \, w . \nabla) \,\,
 \begin{pmatrix} M^{_\ell}  & 0 \\ 0 &  M^{-r}  \end{pmatrix} \, m_d^*(X ,\, t) 
\monendstar 
and the relation (\ref{evolution-moments-decouples}) is demonstrated. \hfill$\square$ 

\bigskip \noindent
How can we define the conserved variables for a bipoint? On the surface,
we have twice as many conserved variables as needed. We return to a basic principle
of lattice Boltzmann schemes: conserved variables are invariant during the relaxation step. 
Since the relaxation operator $ \, m \longmapsto m^* = J_0 \,m \, $ is linear,
it suffices to study the vectors that are invariant under this transformation. 

\bigskip 
{\bf  Proposition 2. Conserved variables for a bipoint}

\noindent
We consider the matrix $ \, J_0 \, $ defined in (\ref{J0-thermique}), its duplicate
\moneq \label{J0-decouple}
J_{0d} = \begin{pmatrix} J_0 & 0 \\ 0 & J_0 \end{pmatrix}
\monend
and the iteration matrix for a null time step
\moneq \label{J0-tilda}
\widetilde{J_0} = \begin{pmatrix} M^r & 0 \\ 0 & M^\ell \end{pmatrix} \, \Pi_\sigma \,\,
 \begin{pmatrix} M^{_\ell}  & 0 \\ 0 &  M^{-r}  \end{pmatrix} \, J_{0d} \, . 
\monend
The eigenvalues of $ \, \widetilde{J_0} \, $ admit the form
%
\moneqstar
+1 \, {\rm (simple)} \,,\,\, 1\!-\!s_j  {\rm (double)} \,,\,\,  1\!-\!s_s  {\rm (simple)} \,,\,\, -1\!+\!s_j  {\rm (double)}  \,,\,\,
 \xi_-  {\rm (simple)} \,,\,\,  \xi_+  {\rm (simple)} 
\monendstar 
where $ \, \xi_{\pm} \, $ are the roots of the equation
\moneqstar
2 \, \xi^2 + s_e \, (\alpha +1) \, \xi + 2\, (s_e-1) = 0 \, . 
\monendstar 
The eigenvector corresponding to the eigenvalue $ \, +1 \, $ is the density
\moneqstar
\rho = \big( 1 ,\, 0 ,\, 0 ,\, \alpha \, \lambda^2  , \, 1 ,\, 0 ,\, 0 ,\, \alpha \, \lambda^2 \big)^{\rm t}
\monendstar 
and the following two eigenvectors associated to the eignevalue $ \, 1\!-\!s_j \, $ define the two components of the momentum: 
\moneqstar \left\{    \begin{array} {rcl}
  j_x &=& \big( 0 ,\, 1 ,\, 0 ,\, 0 ,\, 0 ,\, 1 ,\, 0 ,\, 0 \big)^{\rm t} \\ 
  j_y &=& \big( 0 ,\, 0 ,\, 1 ,\, 0 ,\, 0 ,\, 0 ,\, 1 ,\, 0 \big)^{\rm t} \, . 
  \end{array} \right. \monendstar 

\bigskip \monitem
{\bf  Proof of Proposition 2}

\noindent
After an elementary formal calculus (in ``SageMath'') detailed in \cite{Dubois-zenodo}, we have
%
%
\moneqstar  \widetilde{J_0} = \begin{pmatrix} 
  {1\over4} (1-s_e \, \alpha)  & 0 & 0 &  {{s_e-1}\over{4 \, \lambda^2}}  & 0 & 0 & {{3 + s_e \, \alpha}\over{4}} &  {{1-s_e}\over{4 \, \lambda^2}} \\
  0 & 0 & 0 & 0 & 0 & 1-s_j  & 0 & 0 \\
  0 & 0 & 0 & 0 & 0 & 0 & 1-s_j  & 0 \\ 
  {3\over4} \lambda^2 \, (s_e \, \alpha -1) & 0 & 0 &  {3\over4} (1- s_e) &  {3\over4} \lambda^2 \, (s_e \, \alpha + 3) 
  & 0 & 0 & {{1-s_e}\over{4}}   \\ 
  {{\alpha\, s_e +3}\over{4}} & 0 & 0 &  {{1-s_e}\over{4 \, \lambda^2}} & {{1 - s_e \,\alpha}\over{4}} & 0 & 0 & {{s_e - 1}\over{4 \, \lambda^2}}  \\
 0 &  1-s_j & 0 & 0 & 0 & 0 & 0 & 0 \\ 
 0 & 0 & 1-s_j  & 0 & 0 & 0 & 0 & 0 \\ 
 {1\over4} \lambda^2 \, (s_e \, \alpha +3) & 0 & 0 & {1\over4} (1-s_e) &    {3\over4} \lambda^2 \, (s_e \, \alpha -1)  & 0 & 0 & {3\over4} (1-s_e) 
\end{pmatrix}
\monendstar 

The relations  $ \, \widetilde{J_0}  \, \rho = \rho $ ,  $ \, \widetilde{J_0}  \,\, j_x = (1-s_j) \, j_x \, $
and $ \, \widetilde{J_0}  \,\, j_y = (1-s_j) \, j_y \, $ are elementary. \hfill $\square$ 


\bigskip 
Due to the permutation matrix $ \, \Pi_\sigma $, the structure of the discrete evolution
equation (\ref{evolution-moments-decouples}) is not exactly of the form (\ref{evolution-du22}).
We introduce a ``breathing matrix'' $ \, R_d \, $ in order to consider a conjugate matrix
of the time iteration matrix $ \, \exp (- \Delta t \, w.\nabla) $.

\bigskip 
{\bf  Proposition 3. Breathing matrix}

We define the ``breathing matrix'' $ \, R_d \, $ by the relation
\moneq \label{matrice-respiration} 
\begin{pmatrix} M^r & 0 \\ 0 & M^\ell \end{pmatrix} \, \Pi_\sigma = R_d \,  \begin{pmatrix} M^l & 0 \\ 0 & M^r \end{pmatrix} 
\monend
Then we have
\moneq \label{expression-matrice-respiration}
R_d =  \begin{pmatrix} {1\over4} & 0 & 0 & -{1\over{4 \, \lambda^2}} & {3\over4} & 0 & 0 & {1\over{4 \, \lambda^2}} \\
0 & 0 & 0 & 0 & 0 & 1 & 0 & 0  \\ 
0 & 0 & 0 & 0 & 0 & 0 & 1 & 0  \\
-{3\over4} \lambda^2 & 0 & 0 &  {3\over4} & {3\over4} \lambda^2 & 0 & 0 & {1\over4} \\
{3\over4} & 0 & 0 & {1\over{4 \, \lambda^2}} & {1\over4}  & 0 & 0 & -{1\over{4 \, \lambda^2}} \\
0 & 1 & 0 & 0 & 0 & 0 & 0 & 0 \\ 
0 & 0 & 1 & 0 & 0 & 0 & 0 & 0 \\
{3\over4} \lambda^2 & 0 & 0 &  {1\over4}  & -{3\over4} \, \lambda^2 & 0 & 0 & {3\over4} 
\end{pmatrix} \monend
and the evolution equation admits in the linear case the form
\moneq \label{evolution-moments-decouples-bis}
m_d(X,\, t+\Delta t)  = R_d \,   \begin{pmatrix} M^l & 0 \\ 0 & M^r \end{pmatrix} \,\,
\exp (-\Delta t \,\, w.\nabla) \,\, \begin{pmatrix} M^{-\ell} & 0 \\ 0 & M^{-r}  \end{pmatrix} \, J_{0d} \,\, m_d(X,\, t) \,.
\monend
with $ \, J_{0d}  \, $ defined in (\ref{J0-decouple}).

\bigskip \monitem
{\bf  Proof of Proposition 3}

\noindent
The calculus of the matrix $ \, R_d \, $ is elementary: 
\moneqstar
R_d = \begin{pmatrix} M^r & 0 \\ 0 & M^\ell \end{pmatrix} \,\, \Pi_\sigma  \,\, 
\begin{pmatrix} M^{-\ell} & 0 \\ 0 & M^{-r}  \end{pmatrix} 
\monendstar
and the relation (\ref{expression-matrice-respiration}) is clear.
Moreover, in the linear case, the relations (\ref{relaxation-moments}),
(\ref{J0-thermique}), (\ref{J0-acoustique}), and~(\ref{J0-decouple}) show that
\moneq \label{m-star-decouple} 
m_d^* =  J_{0d} \,\,  m_d \, . 
\monend 
So  the relation (\ref{evolution-moments-decouples-bis}) 
is a direct consequence of  (\ref{evolution-moments-decouples}) and (\ref{expression-matrice-respiration}), and
(\ref{m-star-decouple}). \hfill $\square$ 

\bigskip \noindent 
Next, we eliminate this breathing matrix by calculating its inverse step by step.
We introduce the notations $ \, {\rm I}_k \, $ for the identity matrix of order $ \, k $.
For example,
\moneqstar
{\rm I}_4 = \begin{pmatrix} 1 & 0 & 0 & 0 \\ 0 & 1 & 0 & 0 \\ 0 & 0 & 1 & 0 \\ 0 & 0 & 0 & 1 \end{pmatrix} \, . 
\monendstar
We set 
\moneqstar
N_0 =  \begin{pmatrix} {\rm I}_4 & {\rm I}_4  \\ -{\rm I}_4  & {\rm I}_4   \end{pmatrix} \,, \quad 
N_1 =  \begin{pmatrix} 1 & 0 & 0 & 0 & 0 & 0 & 0 & 0 \\ 0 & 1 & 0 & 0 & 0 & 0 & 0 & 0 \\ 0 & 0 & 1 & 0 & 0 & 0 & 0 & 0 \\ 
0 & 0 & 0 & 1 & 0 & 0 & 0 & 0 \\ 0 & 0 & 0 & 0 & {1\over4} & 0 & 0 & -{1\over{4 \, \lambda^2}} \\  0 & 0 & 0 & 0 & 0 & 1 & 0 & 0 \\ 
0 & 0 & 0 & 0 & 0 & 0 & 1 & 0 \\ 0 & 0 & 0 & 0 & {3\over4} & 0 & 0 & {1\over{4 \, \lambda^2}}  \end{pmatrix} \,.
\monendstar

\smallskip \smallskip \noindent 
Then we have the following  elementary calculus relative to a conjugate of the matrix $\, R_d \,$
evaluated in~(\ref{expression-matrice-respiration}): 
\moneqstar
N_1 \,\, N_0 \,\, R_d \,\, N_0^{-1} \,\, N_1^{-1}   = R_p \equiv {\rm diag} \big( 1 ,\, 1 ,\, 1 ,\, 1 ,\, 1 ,\, -1 -1 ,\, -1 \big) \,. 
\monendstar
We remark that $ \, R_p^2 = {\rm I}_8 \, $ and in consequence $ \, (R_p)^{-1} = R_p $. 
We introduce the ``breathing moments'' $ \, m_r \,$ according to the relation
\moneq \label {moments-respirants}
m_r = (R_p)^{-1} \, N_1 \, N_0 \,\,\,  m_d \, . 
\monend

\bigskip 
{\bf  Proposition 4: discrete evolution of the breathing moments}

With the above notations, the discrete evolution of the breathing moments is given by the relation 
\moneq \label{evolution-moments-respirants}
m_r(X,\, t+\Delta t)  =  \exp \big(-\Delta t \, \Lambda_r \big) \,  J_{0r} \, m_r(X,\, t) \,
\monend
with
\moneq \label{Lambda-r}
 \Lambda_r = N_1 \,\, N_0  \,\, \begin{pmatrix} M^l & 0 \\ 0 & M^r \end{pmatrix} \,\,
\big( w.\nabla \big) \,\, \begin{pmatrix} M^{-\ell} & 0 \\ 0 & M^{-r} \end{pmatrix} \,\,  N_0^{-1} \,\, N_1^{-1} 
\monend
\moneq \label{J0r}
J_{0r} =  N_1 \,\, N_0 \,\,  J_{0d} \,\,   N_0^{-1} \,\, N_1^{-1} \,\, R_p \, .
\monend

\bigskip \monitem
{\bf  Proof of Proposition 4}

\noindent
We have the relation (\ref{evolution-moments-decouples-bis}) and the following calculus

\smallskip \noindent
$ m_r(X,\, t+\Delta t) = (R_p)^{-1} \, N_1 \, N_0 \,\,\,  m_d (X,\, t+\Delta t) $ \qquad due to (\ref{moments-respirants})

\smallskip \noindent \quad
$ = R_p^{-1} \, N_1 \, N_0 \,\,\,  R_d \,   \begin{pmatrix} M^l & 0 \\ 0 & M^r \end{pmatrix} \,\,
\exp (-\Delta t \,\, w.\nabla) \,\, \begin{pmatrix} M^{-\ell} & 0 \\ 0 & M^{-r}  \end{pmatrix} \, J_{0d} \,\, m_d(X,\, t) $ 

\smallskip \noindent \quad
$ = \big(N_1 \, N_0 \, R_d^{-1} \, N_0^{-1} \,  N_1^{-1} \big) \, N_1 \, N_0 \,  R_d \begin{pmatrix} M^l & 0 \\ 0 & M^r \end{pmatrix} 
\exp (-\Delta t \,\, w.\nabla) \begin{pmatrix} M^{-\ell} & 0 \\ 0 & M^{-r}  \end{pmatrix} J_{0d} \,\, m_d(X,\, t) $ 

\hfill due to the definition of $ \, R_p $ 

\smallskip \noindent \quad
$ = \Big[ N_1 \, N_0 \, \begin{pmatrix} M^l & 0 \\ 0 & M^r \end{pmatrix} \exp (-\Delta t \,\, w.\nabla)
\begin{pmatrix} M^{-\ell} & 0 \\ 0 & M^{-r}  \end{pmatrix}   N_0^{-1} \,  N_1^{-1} \Big] \,  N_1 \, N_0 \, J_{0d} \,\, m_d(X,\, t) $ 

\smallskip \noindent \quad
$ = \exp \big( - \Delta t \,  \Lambda_r \big) \,  N_1 \, N_0 \, J_{0d} \,\,  m_d(X,\, t) $ 
\qquad due to (\ref{Lambda-r})

\smallskip \noindent \quad
$ = \exp \big( - \Delta t \,  \Lambda_r \big) \,  N_1 \, N_0 \, J_{0d} \, \big(  N_0^{-1} \,\, N_1^{-1} \, R_p \,\, m_r(X,\, t) \big) $   
\qquad due to (\ref{moments-respirants})

\smallskip \noindent \quad
$ = \exp \big( - \Delta t \,  \Lambda_r \big) \, \big(  N_1 \, N_0 \, J_{0d} \,  N_0^{-1} \,\, N_1^{-1} \, R_p \big) \,\,   m_r(X,\, t) \big) $

\smallskip
and due to (\ref{J0r}), the relation (\ref{evolution-moments-respirants}) is established.
\hfill $\square$ 

\bigskip \bigskip    \noindent {\bf \large    4) \quad  {Equivalent partial differential equations}}

\smallskip \noindent
We first note that if we ignore the permutation operator $ \, \Pi_\sigma  \,$
for the global definition of the D2T4 lattice Boltzmann scheme, 
{\it id est}  if we simplify equations
(\ref{g-droite-j}) and (\ref{g-gauche-j}) by assuming them to be valid for all velocities,
for $ \, 0 \leq j \leq 3 $,  and not only for $ \, 1 \leq j \leq 3 $, 
then we obtain an unsatisfactory result. We refer to Appendix C for such an approach. 

\smallskip \noindent
With the discrete evolution (\ref{evolution-moments-respirants}), we have an equation of the type 
(\ref{evolution-du22}). We can explicit the operator matrix $ \, \Lambda_r \, $ and the matrix matrix
$ \, J_{0r}  $. The matrix  $ \, \Lambda_r \, $ does not depend on the model chosen between diffusion and acoustics.
We have 
\moneqstar 
\Lambda_r  = \begin{pmatrix}
0 & \partial_x & \partial_y & 0 & 0 & 0 & 0 & -\lambda\, \partial_x \\
{3\over8} \lambda^2 \partial_x & 0 & 0 & {1\over8}\, \partial_x & 0 & -{1\over2} \lambda \, \partial_x & -{1\over2} \lambda \, \partial_y & 0 \\ 
{3\over8} \lambda^2 \partial_y & 0 & 0 & {1\over8}\, \partial_y & 0 & -{1\over2} \lambda \, \partial_y & -{3\over2} \lambda \, \partial_x & 0 \\ 
0 & \lambda^2 \partial_x &  \lambda^2 \partial_y & 0 & 0 & 0 & 0 & -\lambda^3 \partial_x \\
0 & 0 & 0 & 0 & 0 & 0 & 0 & 0 \\ 
0 & -{1\over2} \lambda \, \partial_x &  -{1\over2} \lambda \, \partial_y & 0 & 0 & 0 & 0 & {1\over2} \lambda^2 \partial_x \\
0 & -{1\over2} \lambda \, \partial_y &  -{3\over2} \lambda \, \partial_x & 0 & 0 & 0 & 0 & {1\over2} \lambda^2 \partial_y \\
-{3\over4} \lambda \, \partial_x & 0 & 0 & -{1\over{4 \, \lambda}} \, \partial_x & 0 & \partial_x & \partial_y & 0  \end{pmatrix} \, . 
\monendstar 
%
%
We observe that the first row of the matrix  $ \, \Lambda_r \, $ has a non-zero element
in the last position of the first row. This is unusual since in previous calculations (see {\it e.g.} \cite{DL23}),
this matrix element was always zero.
Moreover, due to the respect of the conserved variables of the scheme, the matrix $ \, J_{0r}  \, $ has a structure of the type
\moneq \label{J0r-structure}
J_{0r} = \begin{pmatrix} {\rm I} & 0 \\ S \, E &  {\rm I} - S \end{pmatrix} 
\monend
and the dimensions of this bock decomposition depend on the physical application.

\smallskip  \monitem
For diffusion , we have only one conservation law.
Consequently, the equilibrium matrix~$ \, E_t \, $ is a columns with  7 lines. We have 
\moneq \label{E-thermique}
E_t = \big( 0 ,\, 0 ,\, \lambda^2 \, \alpha ,\, 0 ,\, 0 ,\, 0 ,\, 0  \big)^{\rm t} \, . 
\monend
Moreover, the matrix $ \, S_t \, $ is 7 by 7: 
\moneqstar 
S_t =  \begin{pmatrix} s_j & 0 & 0 & 0 & 0 & 0 & 0 \\ 0 & s_j & 0 & 0 & 0 & 0 & 0 \\  0 & 0 & s_e & 0 & 0 & 0 & 0 \\
0 & 0 & 0 & {1\over4} s_e \, (\alpha+3) & 0 & 0 & {1\over4} s_e \, (1-\alpha) \\ 0 & 0 & 0 & 0 & 2-s_j & 0 & 0 \\
0 & 0 & 0 & 0 & 0 & 2-s_j & 0 \\ 0 & 0 & 0 & -{1\over4} s_e \, (\alpha+3) & 0 & 0 & 2 - {1\over4} s_e \, (1-\alpha) \end{pmatrix} \, . 
\monendstar   
%
%
This matrix does not have the usual structure of a diagonal matrix with elements
that do not depend on equilibrium, as initially proposed by d'Humi\`eres \cite{DDH92}.
The diagonal structure is lost. Moreover, the coefficients
of the matrix $ \, S_t \, $ depend explicitely on the parameter $ \, \alpha $,
which completely determines the choice of equilibrium distribution here,
as shown in (\ref{equilibre-thermique}). Nevertheless,
when the coefficients $ \, \alpha + 3 $, $ \, s_j \, $ and $ \, s_e \, $ are all nonzero, 
this matrix is invertible and we can define the H\'enon \cite{He87} matrix by the relation 
$ \, \Sigma_t \equiv   S_t^{-1} - {1\over2} \, {\rm I}_7  $.
With the notations 
\moneq \label{sigmaj-sigmae}
\sigma_j = {{1}\over{s_j}} - {1\over2}  \,,\,\, \sigma_e = {{1}\over{s_e}} - {1\over2} \, , 
\monend
we have 
\moneq \label{Henon-thermique}
\Sigma_t =  \begin{pmatrix} \sigma_j & 0 & 0 & 0 & 0 & 0 & 0 \\ 0 & \sigma_j & 0 & 0 & 0 & 0 & 0 \\  0 & 0 & \sigma_e & 0 & 0 & 0 & 0 \\
0 & 0 & 0 & {{4\, s_e}\over{\alpha+3}} & 0 & 0 & {{\alpha + 1}\over{2\, (\alpha+3)}} \\ 0 & 0 & 0 & 0 & {{1}\over{4\, \sigma_j}} & 0 & 0 \\
0 & 0 & 0 & 0 & 0 & {{1}\over{4\, \sigma_j}} & 0 \\ 0 & 0 & 0 & {1\over2} & 0 & 0 & 0 \end{pmatrix}  \,. 
\monend   
%
%
We now have all the ingredients to develop the formal analysis ``ABCD'' proposed
by Dubois in \cite{ADGL14,Du22}. We first decompose the matrix $ \, \Lambda_r \, $ into four blocks:
\moneq \label{ABCD}
 \Lambda_r = \begin{pmatrix} A & B \\ C & D \end{pmatrix} \, .
\monend   
For diffusion, $ \, A \, $ is a 1 by 1 matrix and $ \, B \, $ a 1 by 7 one. We have

\vskip -.4 cm 
\moneq \label{A-B-thermique}
A = (0) \,,\,\,
B = \big(  \partial_x ,\, \partial_y ,\, 0 ,\, 0 ,\, 0,\, 0 ,\,  -\lambda\, \partial_x  \big) \, . 
\monend
We have also a 7 by 1 matrix  $\, C \,  $ and a 7 by 7 matrix $ \, D $: 
\moneqstar
C =  \begin{pmatrix} {3\over8} \lambda^2 \partial_x \\ 
{3\over8} \lambda^2 \partial_y  \\  0  \\ 0 \\ 0 \\ 0 \\  -{3\over4} \lambda \, \partial_x \end{pmatrix} \,,\,\,
D =  \begin{pmatrix} 0 & 0 & {1\over8}\, \partial_x & 0 & -{1\over2} \lambda \partial_x & -{1\over2} \lambda \partial_y & 0 \\ 
0 & 0 & {1\over8}\, \partial_y & 0 & -{1\over2} \lambda \partial_y & -{3\over2} \lambda \partial_x & 0 \\ 
\lambda^2 \partial_x &  \lambda^2 \partial_y & 0 & 0 & 0 & 0 & -\lambda^3 \partial_x \\
0 & 0 & 0 & 0 & 0 & 0 & 0 \\ 
-{1\over2} \lambda \,  \partial_x &  -{1\over2} \lambda \,  \partial_y & 0 & 0 & 0 & 0 & {1\over2} \lambda^2 \partial_x \\
-{1\over2} \lambda \,  \partial_y &  -{3\over2} \lambda \,  \partial_x & 0 & 0 & 0 & 0 & {1\over2} \lambda^2 \partial_y \\
0 & 0 & -{1\over{4 \, \lambda}} \, \partial_x  & 0 & \partial_x & \partial_y & 0  \end{pmatrix} \, .
\monendstar

\smallskip  \monitem
In acoustics, we have three conserved variables and five non-conserved moments.
In consequence, the bloc matrices $ \, E\, $  and $ \, S \, $  in (\ref{J0r-structure}) 
have the following dimensions: $ \, 5 \times 3 \, $ and $ \, 5 \times 5 $.
We have
\moneq \label{E-acoustique}
E_a = \begin{pmatrix} \lambda^2 \, \alpha & 0 & 0 \\ 0 & 0 & 0 \\ 0 & 0 & 0 \\  0 & 0 & 0 \\  0 & 0 & 0  \end{pmatrix} 
\monend
%
%
and
\moneqstar 
S_a =  \begin{pmatrix} s_e & 0 & 0 & 0 & 0 \\
0 & {1\over4} s_e \, (\alpha+3) & 0 & 0 & {1\over4} s_e \, (1-\alpha) \\ 0 & 0 & 2 & 0 & 0 \\
0 & 0 & 0 & 2 & 0 \\ 0 & -{1\over4} s_e \, (\alpha+3) & 0 & 0 & 2 - {1\over4} s_e \, (1-\alpha) \end{pmatrix} \, . 
\monendstar   
%
%
If $ \, s_e \, (\alpha+3) \not= 0 $, this matrix is invertible and the associated H\'enon matrix
$ \, \Sigma_a \equiv  S_a^{-1} - {1\over2} \, {\rm I}_5 \, $ is obtained by the relation 
\moneq \label{Henon-acoustique}
\Sigma_a  =  \begin{pmatrix} \sigma_e & 0 & 0 & 0 & 0 \\
0 & {{4\, s_e}\over{\alpha+3}} & 0 & 0 & {{\alpha + 1}\over{2\, (\alpha+3)}} \\ 0 & 0 & 0 & 0 & 0 \\ 0 & 0 & 0 & 0 & 0 \\
0 & {1\over2} & 0 & 0 & 0 \end{pmatrix} \, . 
\monend   
%
%
For acoustics, the ABCD  decomposition  of the matrix $ \, \Lambda_r \, $  into four blocks, as in (\ref{ABCD}), 
contains a~3 by 3 matrix~$ \, A \, $ and a 3 by 5 matrix  $\, B $:  

\vskip -.4 cm 
\moneq \label{A-B-acoustique}
A =  \begin{pmatrix} 0 & \partial_x & \partial_y  \\ {3\over8} \lambda^2 \partial_x & 0 & 0  \\ 
{3\over8} \lambda^2 \partial_y & 0 & 0  \\ \end{pmatrix} \,,\,\, 
B =  \begin{pmatrix} 0 & 0 & 0 & 0 & -\lambda\, \partial_x \\
{1\over8}\, \partial_x & 0 & -{1\over2} \lambda \partial_x & -{1\over2} \lambda \partial_y & 0 \\ 
{1\over8}\, \partial_y & 0 & -{1\over2} \lambda \partial_y & -{3\over2} \lambda \partial_x & 0 \end{pmatrix} \, .  
\monend
We have also a 5 by 3 matrix  $\, C \,  $ and a 5 by 5 matrix $ \, D $: 
\moneqstar
C =  \begin{pmatrix} 
0 & \lambda^2 \partial_x &  \lambda^2 \partial_y \\ 0 & 0 & 0  \\ 0 & -{1\over2} \lambda \, \partial_x &  -{1\over2} \lambda \, \partial_y  \\
0 & -{1\over2} \lambda \, \partial_y &  -{3\over2} \lambda \,  \partial_x  \\ -{3\over4} \lambda \, \partial_x & 0 & 0   \end{pmatrix} \,,\,\,
D =  \begin{pmatrix}  0 & 0 & 0 & 0 & -\lambda^3 \partial_x \\ 0 & 0 & 0 & 0 & 0 \\ 0 & 0 & 0 & 0 & {1\over2} \lambda^2 \partial_x \\
0 & 0 & 0 & 0 & {1\over2} \lambda^2 \partial_y \\  -{1\over{4 \, \lambda}} \, \partial_x  & 0 & \partial_x & \partial_y & 0  \end{pmatrix} \, . 
\monendstar 
%
%

\bigskip 
{\bf  Proposition 5. Second order equivalent partial differential equation for diffusion}

In the diffusive case with only one conservation moment and equilibrium given by the relation
(\ref{equilibre-thermique}), the equivalent partial differential equation can be written
\moneq \label{edp-thermique}
{{\partial \rho}\over{\partial t}} - {{\lambda^2}\over{8}} \, (\alpha+3) \,\Delta t \, \sigma_j \,
\big( \partial_x^2 +  \partial_y^2 \big) \rho = {\rm O}(\Delta t^2) \,. 
\monend   
%
%

\newpage 
\bigskip \monitem
{\bf  Proof of Proposition 5}

\noindent
With all the hypotheses recalled previously, it is know that in the linear case, the partial differential
equation at order 2 takes the form 
\moneq \label{equivalent-edp-ordre-2}
\partial_{t} W + \alpha_1 \, W + \Delta t \, \alpha_2 \, W =  {\rm O} (\Delta t^2)
\monend
with operators $ \, \alpha_1 \, $ and $ \, \alpha_2 \, $ obtained by the algorithm \cite{ADGL14,Du22}: 
\moneq \label{algorithme-berlin} \left\{ \begin{array}   {l}
\alpha_1 = A + B \, E  \\
\beta_1 =  E \, \alpha_1 - (C +  D \, E)   \\ 
\alpha_2 =  B \,  \Sigma \, \beta_1   \, . 
\end{array} \right. \monend
In the diffusion case, $ \, W = \rho \, $ is a scalar field and is is clear from (\ref{A-B-thermique}) 
that 
\moneqstar 
\alpha_1 = A + B \, E_t = 0 \, . 
\monendstar
At second order, we have

\smallskip \noindent
$ C + D \, E_t = C  + 
\begin{pmatrix} 0 & 0 & {1\over8}\, \partial_x & 0 & -{1\over2} \lambda \partial_x & -{1\over2} \lambda \partial_y & 0 \\ 
0 & 0 & {1\over8}\, \partial_y & 0 & -{1\over2} \lambda \partial_y & -{3\over2} \lambda \partial_x & 0 \\ 
\lambda^2 \partial_x &  \lambda^2 \partial_y & 0 & 0 & 0 & 0 & -\lambda^3 \partial_x \\
0 & 0 & 0 & 0 & 0 & 0 & 0 \\ 
-{1\over2} \lambda \,  \partial_x &  -{1\over2} \lambda \,  \partial_y & 0 & 0 & 0 & 0 & {1\over2} \lambda^2 \partial_x \\
-{1\over2} \lambda \,  \partial_y &  -{3\over2} \lambda \,  \partial_x & 0 & 0 & 0 & 0 & {1\over2} \lambda^2 \partial_y \\
0 & 0 & -{1\over{4 \, \lambda}} \, \partial_x  & 0 & \partial_x & \partial_y & 0  \end{pmatrix}
\begin{pmatrix} 0 \\ 0 \\ \lambda^2 \, \alpha  \\  0 \\  0 \\  0 \\  0  \end{pmatrix} $

\smallskip \noindent \qquad \qquad 
$ \,\,\, = \begin{pmatrix} {3\over8} \lambda^2 \partial_x \\ 
{3\over8} \lambda^2 \partial_y  \\  0  \\ 0 \\ 0 \\ 0 \\  -{3\over4} \lambda \, \partial_x \end{pmatrix}  + 
\begin{pmatrix} {1\over8} \lambda^2 \, \alpha \, \partial_x \\  {1\over8} \lambda^2 \, \alpha \, \partial_y \\
  0 \\ 0 \\ 0 \\ 0 \\ 0 \end{pmatrix} =
\begin{pmatrix} {1\over8} \lambda^2 \, (\alpha+3) \, \partial_x \\ {1\over8} \lambda^2 \, (\alpha+3) \, \partial_y
 \\  0  \\ 0 \\ 0 \\ 0 \\  -{3\over4} \lambda \, \partial_x \end{pmatrix}  $ 

\smallskip \smallskip \noindent
and the operator $ \, \beta_1 \, $ takes the expression 
\moneqstar
\beta_1 = \big( -{1\over8} \, \lambda^2 \, (\alpha+3) \, \partial_x \,,\,\,   -{1\over8} \, \lambda^2 \, (\alpha+3) \, \partial_y  \,,\,\, 
  0 \,,\,\,  0 \,,\,\, 0 \,,\,\, 0 \,,\,\, {1\over4} \, \lambda  \, (\alpha+3) \, \partial_x \big)^{\rm t} \, .
\monendstar 
%
%
With the operator $ \, B \, $ given in (\ref{A-B-thermique})  and the H\'enon matrix $ \, \Sigma \, $
in (\ref{Henon-thermique}), the expression of $ \, \alpha_2 \, $ is computed as follows:

\smallskip \noindent
$  \alpha_2 = B \, \Sigma_t \, \beta_1 = \Big( \partial_x  \,\,\,\, \partial_y  \,\,\,\, 0 \,\,\,\, 0 \,\,\,\, 0 \,\,\,\, 0
\,\,\,\, -\lambda\, \partial_x \Big)
\begin{pmatrix} \sigma_j & 0 & 0 & 0 & 0 & 0 & 0 \\ 0 & \sigma_j & 0 & 0 & 0 & 0 & 0 \\  0 & 0 & \sigma_e & 0 & 0 & 0 & 0 \\
0 & 0 & 0 & {{4\, s_e}\over{\alpha+3}} & 0 & 0 & {{\alpha + 1}\over{2\, (\alpha+3)}} \\ 0 & 0 & 0 & 0 & {{1}\over{4\, \sigma_j}} & 0 & 0 \\
0 & 0 & 0 & 0 & 0 & {{1}\over{4\, \sigma_j}} & 0 \\ 0 & 0 & 0 & {1\over2} & 0 & 0 & 0 \end{pmatrix} 
\beta_1 $

\smallskip \noindent \qquad
$ = \Big(   \sigma_j \,  \partial_x  \,\,\,\,  \sigma_j \,  \partial_y   \,\,\,\, 0 \,\, 0 \,\, 0 \,\,\,\, 0
\,\,\,\, 0 \Big) 
\begin{pmatrix}  -{1\over8} \, \lambda^2 \, (\alpha+3) \, \partial_x  \\   -{1\over8} \, \lambda^2 \, (\alpha+3) \, \partial_y   \\  
  0 \\  0 \\ 0 \\ 0 \\  {1\over4} \, \lambda  \, (\alpha+3) \, \partial_x \end{pmatrix}  =
 - {{1}\over{8}} \, \lambda^2 \, (\alpha+3) \, \sigma_j \, \big( \partial_x^2 +  \partial_y^2 \big) $ 

\smallskip \noindent 
and the relation (\ref{edp-thermique}) is established.
The proposition 5 is proven. \hfill $\square$ 

\bigskip 
{\bf  Proposition 6. Second order equivalent partial differential equation for acoustics}

In the acoustic case, we have three conserved scalar quantities, {\it id est} the scalar $ \,\rho \, $
and the momentum $ \, J = \big( j_x ,\, j_y \big) $. The equilibrium of the last moment
satisfies the relation (\ref{equilibre-acoustique}).
The system of  equivalent partial differential equations  satisfies the relations 
\moneq \label{edp-acoustique}  \left\{    \begin{array} {rcl}
 {{\partial \rho}\over{\partial t}} + {\rm div} J =  {\rm O}(\Delta t^2) \\  
 {{\partial J}\over{\partial t}} + c_0^2 \, \nabla \rho - \zeta \, \nabla ( {\rm div} J ) = {\rm O}(\Delta t^2)
\end{array} \right. \monend
with $ \,  {\rm div} J \equiv \partial_x j_x + \partial_y j_y \, $ and 
\moneq \label{c-son-zeta} 
c_0^2 =  {{\lambda^2}\over{8}} \, (\alpha+3) \,,\,\,\,\,
\zeta  = {{\lambda^2}\over{8}} \, \Delta t \, (1-\alpha)\, \sigma_e \,,\,\,\,\,
\sigma_e = {1\over{s_e}} - {1\over2} \, . 
\monend
The viscous term in the momentum equation in (\ref{edp-acoustique}) corresponds to a bulk viscosity; 
the shear viscosity of this model is equal  to zero.

\bigskip \noindent 
The defects noted in Appendices B and C, with an overly superficial analysis, are no longer present. 

\bigskip \monitem
{\bf  Proof of Proposition 6}

\noindent
From the relations (\ref{algorithme-berlin})

\smallskip
$ \alpha_1 = A + B \, E_a =
\begin{pmatrix} 0 & \partial_x & \partial_y  \\ {3\over8} \lambda^2 \partial_x & 0 & 0  \\ 
{3\over8} \lambda^2 \partial_y & 0 & 0 \\ \end{pmatrix} + 
\begin{pmatrix} 0 & 0 & 0 & 0 & -\lambda\, \partial_x \\
{1\over8}\, \partial_x & 0 & -{1\over2} \lambda \partial_x & -{1\over2} \lambda \partial_y & 0 \\ 
{1\over8}\, \partial_y & 0 & -{1\over2} \lambda \partial_y & -{3\over2} \lambda \partial_x & 0 \end{pmatrix}
 \begin{pmatrix} \lambda^2 \, \alpha & 0 & 0 \\ 0 & 0 & 0 \\ 0 & 0 & 0 \\  0 & 0 & 0 \\  0 & 0 & 0  \end{pmatrix} $

 \qquad \qquad \qquad
 $ \,\,\,\, = \begin{pmatrix} 0 & \partial_x & \partial_y  \\ {1\over8} \lambda^2 \, (\alpha + 3) \, \partial_x & 0 & 0  \\ 
  {1\over8} \lambda^2 \, (\alpha + 3) \, \partial_y & 0 & 0  \\ \end{pmatrix} $.

\smallskip \noindent
Then the first order terms of the acoustics equations (\ref{edp-acoustique}) are justified,
and we have moreover $ \, c_0^2 =  {{\lambda^2}\over{8}} \, (\alpha+3) $. We have on the other hand

\smallskip \noindent 
$ \beta_1 = E_a \, \alpha_1 - (C +  D \, E) $

\smallskip \noindent \quad 
$=  \begin{pmatrix} \lambda^2 \, \alpha & 0 & 0 \\ 0 & 0 & 0 \\ 0 & 0 & 0 \\  0 & 0 & 0 \\  0 & 0 & 0  \end{pmatrix}
\begin{pmatrix} 0 & \partial_x & \partial_y  \\ {1\over8} \lambda^2 \, (\alpha + 3) \, \partial_x & 0 & 0  \\ 
  {1\over8} \lambda^2 \, (\alpha + 3) \, \partial_y & 0 & 0  \\ \end{pmatrix}
- \begin{pmatrix} 
0 & \lambda^2 \partial_x &  \lambda^2 \partial_y \\ 0 & 0 & 0  \\ 0 & -{1\over2} \lambda \, \partial_x &  -{1\over2} \lambda \, \partial_y  \\
0 & -{1\over2} \lambda \, \partial_y &  -{3\over2} \lambda \, \partial_x  \\ -{3\over4} \lambda  \partial_x & 0 & 0   \end{pmatrix} $

\smallskip \noindent \qquad \qquad \qquad 
$ -  \begin{pmatrix}  0 & 0 & 0 & 0 & -\lambda^3 \partial_x \\ 0 & 0 & 0 & 0 & 0 \\ 0 & 0 & 0 & 0 & {1\over2} \lambda^2 \partial_x \\
0 & 0 & 0 & 0 & {1\over2} \lambda^2 \partial_y \\  -{1\over{4 \, \lambda}} \, \partial_x  & 0 & \partial_x & \partial_y & 0  \end{pmatrix} 
\begin{pmatrix} \lambda^2 \, \alpha & 0 & 0 \\ 0 & 0 & 0 \\ 0 & 0 & 0 \\  0 & 0 & 0 \\  0 & 0 & 0  \end{pmatrix} $

\smallskip \noindent \quad 
$ = \begin{pmatrix} 0 & \lambda^2 \, \alpha \, \partial_x &  \lambda^2 \, \alpha \,  \partial_y  \\
  0 & 0 & 0 \\ 0 & 0 & 0 \\  0 & 0 & 0 \\  0 & 0 & 0 \end{pmatrix}  
 - \begin{pmatrix} 
0 & \lambda^2 \partial_x &  \lambda^2 \partial_y \\ 0 & 0 & 0  \\ 0 & -{1\over2} \lambda \, \partial_x &  -{1\over2} \lambda \, \partial_y  \\
0 & -{1\over2} \lambda \, \partial_y &  -{3\over2} \lambda \,  \partial_x  \\ -{3\over4} \lambda \, \partial_x & 0 & 0   \end{pmatrix}
-  \begin{pmatrix}  0 & 0 & 0 \\ 0 & 0 & 0 \\ 0 & 0 & 0 \\ 0 & 0 & 0 \\  -{1\over4} \lambda \, \alpha \,  \partial_x & 0 & 0\end{pmatrix} $

\smallskip \noindent \quad 
$ = \begin{pmatrix} 0 & \lambda^2 \, (\alpha-1) \, \partial_x &  \lambda^2  \, (\alpha-1) \, \partial_y \\ 0 & 0 & 0  \\
0 & {1\over2} \lambda \, \partial_x &  {1\over2} \lambda \, \partial_y  \\
0 & {1\over2} \lambda \, \partial_y &  {3\over2} \lambda \, \partial_x  \\ {3\over4} \lambda \, (\alpha+3) \,  \partial_x & 0 & 0   \end{pmatrix} \, .$
%

\noindent Then

\noindent
$ \alpha_2 = B \, \Sigma_a \, \beta_1 =
\begin{pmatrix} 0 & 0 & 0 & 0 & -\lambda\, \partial_x \\
{1\over8}\, \partial_x & 0 & -{1\over2} \lambda \partial_x & -{1\over2} \lambda \partial_y & 0 \\ 
{1\over8}\, \partial_y & 0 & -{1\over2} \lambda \partial_y & -{3\over2} \lambda \partial_x & 0 \end{pmatrix}
\begin{pmatrix} \sigma_e & 0 & 0 & 0 & 0 \\
0 & {{4\, s_e}\over{\alpha+3}} & 0 & 0 & {{\alpha + 1}\over{2\, (\alpha+3)}} \\ 0 & 0 & 0 & 0 & 0 \\ 0 & 0 & 0 & 0 & 0 \\
0 & {1\over2} & 0 & 0 & 0 \end{pmatrix} \, \beta_1 $

\smallskip \noindent \qquad 
$ =  \begin{pmatrix} 0 & 0 & 0 & 0 & 0 \\
{1\over8}\, \sigma_e \, \partial_x & 0 & 0 & 0 & 0 \\ 
{1\over8}\, \sigma_e \, \partial_y & 0 & 0 & 0 & 0  \end{pmatrix}
 \begin{pmatrix} 0 & \lambda^2 \, (\alpha-1) \, \partial_x &  \lambda^2  \, (\alpha-1) \, \partial_y \\ 0 & 0 & 0  \\
0 & {1\over2} \lambda \, \partial_x &  {1\over2} \lambda \, \partial_y  \\
0 & {1\over2} \lambda \, \partial_y &  {3\over2} \lambda \, \partial_x  \\ {3\over4} \lambda \, (\alpha+3) \,  \partial_x & 0 & 0   \end{pmatrix} $

\smallskip \noindent \qquad 
 $ =  \begin{pmatrix} 0 & 0 & 0  \\
0 & {1\over8}\lambda^2 \, (\alpha-1)\, \sigma_e \, \partial_x^2 &  {1\over8}\lambda^2 \, (\alpha-1)\, \sigma_e \, \partial_x \, \partial_y \\   
0 &  {1\over8}\lambda^2 \, (\alpha-1)\, \sigma_e \, \partial_x \, \partial_y &  {1\over8}\lambda^2 \, (\alpha-1)\, \sigma_e \, \partial_y^2 
\end{pmatrix} =
- {1\over8}\lambda^2 \,(1-\alpha) \,  \begin{pmatrix} 0 & 0 & 0  \\ 0 & \partial_x^2 & \partial_x \, \partial_y \\
0 & \partial_x \, \partial_y &  \partial_y^2  \end{pmatrix} $. 

\smallskip \noindent
The structure of the second order operators in (\ref{edp-acoustique}) is explained.
Moreover, the relation $ \, \zeta  = {{\lambda^2}\over{8}} \, \Delta t \, (1-\alpha) \, \sigma_e $
is also a consequence of the previous calculus. The proposition 6 is proven. \hfill $\square$

\bigskip \bigskip    \noindent {\bf \large    5) \quad  {Numerical experiments}}

\smallskip \noindent
We performed fundamental  numerical experiments in a rectangle with periodic boundary conditions
to avoid any contamination of the results by effects due to the boarders.

\bigskip \monitem Mesh generation

\begin{figure}    [H]  \centering
\centerline{\includegraphics[width=.77\textwidth]{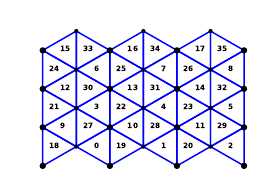}}
\vskip -1.5 cm 
\caption{Mesh with 36 equilateral triangles and associated vertices for $ \, n_x = 3 $.}
\label{fig-d2t4} \end{figure}

In order to be able to mesh the domain with equilateral triangles, we chose
\moneqstar 
\Omega  =  (0,\,\sqrt{3}) \times  (0,\, 1)  \, .
\monendstar
The grids are simply set by an integer $\, n_x$. We place $\, 2 \, n_x \, $  triangles along the $x$-axis and $\, n_x \, $ 
along the $y$-axis, as shown in Figure \ref{fig-d2t4}. 
The total number of triangles is simply $ \, 4 \, n_x^2  $.
For $ \, n_x = 10 $, $\, 20 \, $ and $ \, 40 $, we obtain $ \, 400 $, $ \, 1600 \, $
and $ \, 6400 \, $ triangles respectively. 
The graphics outputs use a relatively coarse rectangular grid.
They are represented with thicker dots in  Figure~\ref{fig-d2t4}.

\smallskip
To implement the D2T4 scheme, it is essential to have a list of neighboring triangles.
This is obtained by following the speed numbering shown in the figure \ref{fig-bi-triangle}.
For example, the four neighbors of triangle 3 are, in this order, the triangles 3, 22, 30 and 27.
To take into account the periodic boundary conditions, 
the neighbors of triangle 16  are 16, 34, 19, 25, and the neighbors of triangle 2 are 2, 18, 29 and 35. 

\vskip -1.5 cm 
\bigskip 

\bigskip \monitem Diffusion   

With the periodic initial density field
\moneqstar 
\rho (x,\, y ,\, 0) = \overline{\rho}(0)  \,\, \sin \big( k_x \, x \big) \, \sin  \big( k_y \, y \big) 
\monendstar 
the solution of the heat equation
\moneqstar 
\partial_t \rho - \kappa \, \big( \partial_x^2 + \partial_y^2 \big)  \rho = 0 
\monendstar 
with 
\moneq \label{kappa-thermique} 
\kappa =  {{\lambda}\over{2}} \, (\alpha+3) \,\Delta t \,\,  \sigma_j
=  {{\lambda}\over{2}} \, (\alpha+3) \,\Delta t \,\, \Big( {1\over{s_j}} - {1\over2} \Big) 
\monend
due to the equation (\ref{edp-thermique}),
is simply an exponential decay:
\moneqstar 
\rho  (x,\, y ,\, 0)  =  \exp (-\kappa \, |k|^2 \, t) \,\, \rho  (x,\, y ,\, 0) 
\monendstar 
with
\moneqstar
|k|^2 = k_x^2 + k_y^2 \, . 
\monendstar

\begin{figure}    [H]  \centering
\centerline{\includegraphics[width=.77\textwidth]{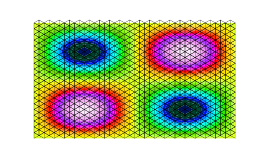}}
\vskip -.9 cm 
\caption{Diffusion test case, $\, \kappa = 10^{-3}$. Density field at final time, 1600  triangles and 20 time steps.
The color map simply allows you to identify the contour lines.
This remark also applies to the following figures.}
\label{fig-thermique-densite} \end{figure}

It is possible to compare this exact solution and the one computed with the lattice Boltzmann D2T4 scheme.
We have chosen 
\moneq \label{parametres-d2t4} 
\alpha = -1 \,,\,\, s_e  = 1.9 \,,\,\, \lambda = 1 \,,\,\, k_x = {{2 \, \pi}\over{\sqrt{3}}}  \,,\,\, k_y = 2 \, \pi \,. 
\monend
The parameter $ \, s_j \, $  is adapted  to the mesh size and to the given value of the diffusivity
$ \, \kappa \, $ through the relation  (\ref{kappa-thermique}). 
A typical result is proposed in Figure \ref{fig-thermique-densite}. The quantitative results for three
nested meshes are displayed in Table \ref{tab-erreurs-thermique}. 

\vskip -.7 cm 

\begin{table}[H]
\smallskip
\centerline {\begin{tabular}{|c|c|c|c|c|c|c|c|c|c|}    \hline 
diffusivity $ \, \kappa $  & $10^{-3}$ & $10^{-3}$ &  $10^{-4}$ & $10^{-4}$ &  $10^{-5}$ & $10^{-5}$   \\   \hline
mesh points  & $ s_j $ & $ \ell^{\infty} \, $ error  & $ s_j $ & $ \ell^{\infty} \, $ error & $ s_j $ & $ \ell^{\infty} \, $ error \\ \hline
 $\, n_x = 10 $, $ \, \Delta x = 0.0577 $   & 1.76  &  1.26 $ 10^{-2}$ &  1.97 &   1.32 $ 10^{-2}$ &  1.997  &   1.29 $ 10^{-2}$ \\  \hline
 $\, n_x = 20 $, $ \, \Delta x = 0.0288 $   & 1.57  &  1.84 $ 10^{-3}$ &  1.95 &   2.88 $ 10^{-3}$  & 1.994  &   2.52 $ 10^{-3}$ \\  \hline
 $\, n_x = 20 $, $ \, \Delta x = 0.0144 $   & 1.29  &  4.82 $ 10^{-4}$ & 1.89 &    6.89 $ 10^{-4}$  & 1.989  &   5.37 $ 10^{-4}$  \\  \hline
convergence order & & 2.35 &&  2.13  && 2.29 \\  \hline
\end{tabular}}
\caption{Diffusion test case. Errors as the mesh size tends to zero and order of convergence.}
\label{tab-erreurs-thermique} \end{table}

\smallskip \noindent
We observe very good quality convergence for these numerical experiments. 
The results for diffusion performance confirm our previous work \cite{DL13}.
If we continue the development to the next order, we can improve the result. The rest in the equivalent equation (\ref{edp-thermique})
is of order~3:
$ \, {{\partial \rho}\over{\partial t}} - {{\lambda}\over{2}} \, (\alpha+3) \,\Delta t \, \sigma_j \, \big( \partial_x^2
+  \partial_y^2 \big) \rho = {\rm O}(\Delta t^3) $.   The proof is in evidence in our ``Sagemath'' software (see \cite{Dubois-zenodo}) 
but is not detailed here. 

\bigskip \monitem Acoustics   

\noindent 
The implementation of a real periodic analytical solution for the system of equations (\ref{edp-acoustique})
requires an algebraic calculation, which is detailed in the following proposition.  

\newpage 
\bigskip 
{\bf  Proposition 7. Exact periodic solution of linear acoustic in a rectangle}

With the initial contion
\moneqstar \left\{    \begin{array} {rcl}
\rho(x,\, y ,\, 0) &=&  \overline{\rho} \, \, \cos(k_x  \, x) \,  \cos(k_y  \, y) \\ 
J_x (x,\, y ,\, 0) &=&  \overline{\rho} \, \, {{k_x}\over{ |k|^2}}  \,  \sin(k_x  \, x) \,  \cos(k_y  \, y) \\ 
J_y (x,\, y ,\, 0) &=&  \overline{\rho} \, \, {{k_y}\over{ |k|^2}}  \,  \cos(k_x  \, x) \,  \sin(k_y  \, y)
\end{array} \right. \monendstar
the solution of the acoustic model (\ref{edp-acoustique}) can the written 
\moneq \label{rho-jj-acoustique}
\left\{    \begin{array} {rcl}
\rho(x,\, y ,\, t) &=&  \overline{\rho} \,\,  \cos(\omega \, t) \, \exp(- \theta \, t) \,  \cos(k_x  \, x) \,  \cos(k_y  \, y) \\ 
J_x (x,\, y ,\, t) &=&  \overline{\rho} \,\,  {{k_x}\over{ |k|^2}}  \,
\big[ \theta \, \cos(\omega \, t) + \omega \, \sin(\omega \, t) \big]  \, \exp(- \theta \, t) \,  \sin(k_x  \, x) \,  \cos(k_y  \, y)  \\ 
J_y (x,\, y ,\, t) &=&  \overline{\rho} \,\,  {{k_y}\over{ |k|^2}}  \,
\big[ \theta \, \cos(\omega \, t) + \omega \, \sin(\omega \, t) \big] \, \exp(- \theta \, t) \,  \cos(k_x  \, x) \,  \sin(k_y  \, y) 
\end{array} \right. \monend
with the parameters  $ \, \omega \, $ and  $ \, \theta $   satisfying the conditions 
\moneq \label{omega-theta-acoustique}
\zeta  \, |k|^2 = 2 \, \theta   \,,\,\, c_0^2  \, |k|^2 = \omega^2 + \theta^2  \,,\,\, \omega > 0  \,,\,\, |k|^2 = k_x^2 + k_y^2  \,.
\monend
We observe that such a solution exists only when 
\moneqstar
c_0 > {1\over2} \, \zeta \, |k| \, .
\monendstar 

\bigskip \monitem
{\bf  Proof of Proposition 7}

\noindent
From the second and third relations of (\ref{rho-jj-acoustique}), we have 
\moneqstar
{\rm div} J =  \overline{\rho} \,\,  \big[ \theta \, \cos(\omega \, t) + \omega \, \sin(\omega \, t) \big]  \, \exp(- \theta \, t)
\,  \cos(k_x  \, x) \,  \cos(k_y  \, y)
\monendstar
and the first equation of (\ref{edp-acoustique})
is a consequence of the identity
\moneqstar
{{\dd}\over{\dd t}} \big(  \cos(\omega \, t) \, \exp(- \theta \, t) \big)
+     \big[ \theta \, \cos(\omega \, t) + \omega \, \sin(\omega \, t) \big]  \, \exp(- \theta \, t) = 0 \, . 
\monendstar
We have also

\smallskip \noindent
$ \partial_t J_x + c_0^2 \, \partial_x J_x - \zeta \,\, \partial_x \big( {\rm div} J \big) $

\smallskip \noindent \qquad 
$ = \overline{\rho} \,   {{k_x}\over{ |k|^2}}  \,  \cos(k_x  \, x) \,  \cos(k_y  \, y) \,\exp(- \theta \, t) \big) \,  \big[
\big( \omega^2 - \theta^2 \big)  \, \cos(\omega \, t) - 2 \, \theta \, \omega \,  \sin(\omega \, t) \big]  \Big] $

\noindent \qquad \qquad \qquad \qquad \qquad \qquad \qquad
$ - \, c_0^2 \, |k|^2 \,  \cos(\omega \, t) + \zeta \,  |k|^2 \,  \big[ \theta \, \cos(\omega \, t) + \omega \, \sin(\omega \, t) \big] $

\smallskip \noindent \qquad 
$ = \overline{\rho} \,   {{k_x}\over{ |k|^2}}  \,  \cos(k_x  \, x) \,  \cos(k_y  \, y) \,\exp(- \theta \, t) \big) \,  \big[
 \big(  \omega^2 - \theta^2   - c_0^2 \, |k|^2 +  \zeta \,  |k|^2 \,  \theta \big) \,  \cos(\omega \, t)    $

\noindent \qquad \qquad \qquad \qquad \qquad \qquad \qquad
$ + \, \omega \, \big(  \zeta \,  |k|^2 - 2 \, \theta \big) \,   \sin (\omega \, t) \Big]  $.

\smallskip \noindent
Due to the first condition of (\ref{omega-theta-acoustique}), the coefficient of $ \,  \sin (\omega \, t) \, $
is zero. Moreover, 
\moneqstar
\omega^2 - \theta^2   - c_0^2 \, |k|^2 +  \zeta \,  |k|^2 \,  \theta =
\omega^2 - \theta^2   - c_0^2 \, |k|^2 + 2 \, \theta^2 = \omega^2 + \theta^2   - c_0^2 \, |k|^2 = 0
\monendstar
due to the second relation of  (\ref{omega-theta-acoustique}).
So the second equation of the system (\ref{edp-acoustique}) is satisfied.
The proof is similar of the third equation of (\ref{edp-acoustique}). 
The proof is completed.
\hfill $\square$

\bigskip 
For the acoustic test case, the parameter $ \, s_j \, $ does not exist anymore.
The bulk viscosity $ \, \zeta \, $ is given according to relation (\ref{c-son-zeta}).
We have fixed three values for this physical parameter:
$ \, \zeta = 10^{-3} $, $ \, 10^{-4} \, $ and $ \, 10^{-5} $.
The parameters $ \, \lambda  $, $  \alpha \, $ and the wave vector $ \, (k_x,\, k_y) \, $
still follow the relation~(\ref{parametres-d2t4}). 
We observe that  $ \, c_0^2 = {1\over4} $.
The parameter $ \, s_e \, $ for the relaxation of the nonconserved moment 
is fixed as a function of $ \, \zeta \, $ and the number
of mesh points in order to satisfy~(\ref{c-son-zeta}).
%

\vskip -.8 cm 
\begin{figure}    [H]  \centering
\centerline{{\includegraphics[width=.62\textwidth]{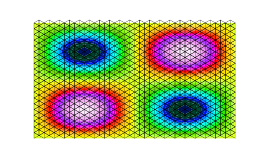}}}
\vskip -.9 cm 
\caption{Acoustics test case, $\, \zeta = 10^{-3}$; density field for 1600  triangles and 20 time steps.}
\label{fig-acoustique-densite-nx20} \end{figure}

\vskip -.8 cm 
\begin{figure}    [H]  \centering
\centerline{{\includegraphics[width=.62\textwidth]{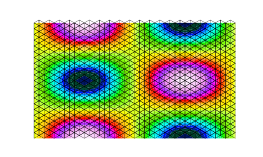}}
\hskip -2.2 cm 
{\includegraphics[width=.62\textwidth]{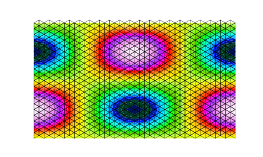}}}
\vskip -.9 cm 
\caption{Acoustics test case, $\, \zeta = 10^{-3}$; first  component of the momentum on the left,
 second  component on the right for  for 1600  triangles and 20 time steps.}
\label{fig-acoustique-jx-jy-nx20} \end{figure}
\bigskip 

\smallskip \noindent
A qualitative view of the results is provided in the figures \ref{fig-acoustique-densite-nx20}
and \ref{fig-acoustique-jx-jy-nx20}.
In the tables \ref{tab-erreurs-acoustique-10-3},  \ref{tab-erreurs-acoustique-10-4} and~\ref{tab-erreurs-acoustique-10-5}, 
we have recalled the values of the $ \, s_e \, $  parameter and explained the relative errors for the density
and the two components of the momentum. The density has good second-order convergence properties, as expected.
Furthermore, the pattern exhibits convergent behavior for both components of the impulsion.
However, the measured convergence is only second order, which leaves a new question open. 

\smallskip 
\begin{table}[H] 
\centerline {\begin{tabular}{|c|c|c|c|c|c|}    \hline 
mesh points  &  $\, s_e $ & $ \rho$  &  $ J_x $ &  $ J_y $  \\   \hline
 $\, n_x = 10 $, $ \, \Delta x = 0.0577 $   &  1.757 & 1.330  $10^{-2}$ &   1.049 $10^{-1}$  &  1.049 $10^{-1}$  \\  \hline
 $\, n_x = 20 $, $ \, \Delta x = 0.0288 $   &  1.566 & 2.974  $10^{-3}$ &   5.230  $10^{-2}$ &  5.230 $10^{-2}$  \\  \hline
 $\, n_x = 20 $, $ \, \Delta x = 0.0144 $   &  1.288 & 4.523  $10^{-4}$ &   2.602 $10^{-2}$  &  2.602 $10^{-2}$  \\  \hline
convergence order & &  2.44   &  1.01  &   1.01   \\  \hline
\end{tabular}}
\caption{Acoustic test case $\, \zeta = 10^{-3}$. Errors for density and momentum as the mesh size tends to zero.}
\label{tab-erreurs-acoustique-10-3} \end{table} 

\begin{table}[H]
\centerline {\begin{tabular}{|c|c|c|c|c|c|}    \hline 
mesh points  &  $\, s_e $ & $ \rho$  &  $ J_x $ &  $ J_y $  \\   \hline
 $\, n_x = 10 $, $ \, \Delta x = 0.0577 $   &  1.973 &  1.326  $10^{-2}$ &   1.055  $10^{-1}$  &    1.055  $10^{-1}$  \\  \hline
 $\, n_x = 20 $, $ \, \Delta x = 0.0288 $   &  1.946 &  3.370  $10^{-3}$ &   5.265  $10^{-2}$  &    5.265  $10^{-2}$  \\  \hline
 $\, n_x = 20 $, $ \, \Delta x = 0.0144 $   &  1.895 &  8.275  $10^{-4}$ &   2.620  $10^{-2}$  &    2.620  $10^{-2}$  \\  \hline
convergence order & &   2.00    &   1.00  &   1.00   \\  \hline
\end{tabular}}
\caption{Acoustic test case $\, \zeta = 10^{-4}$. Errors for density and momentum as the mesh size tends to zero.}
\label{tab-erreurs-acoustique-10-4} \end{table} 

\begin{table}[H]
\centerline {\begin{tabular}{|c|c|c|c|c|c|}    \hline 
mesh points  &  $\, s_e $ & $ \rho$  &  $ J_x $ &  $ J_y $  \\   \hline
 $\, n_x = 10 $, $ \, \Delta x = 0.0577 $   &   1.997  &    1.324  $10^{-2}$ &    1.056  $10^{-1}$  &      1.056 $10^{-1}$  \\  \hline
 $\, n_x = 20 $, $ \, \Delta x = 0.0288 $   &   1.994  &    3.389  $10^{-3}$ &    5.269  $10^{-2}$  &      5.269 $10^{-2}$  \\  \hline
 $\, n_x = 20 $, $ \, \Delta x = 0.0144 $   &   1.989  &    8.330  $10^{-4}$ &    2.622  $10^{-2}$  &      2.622 $10^{-2}$  \\  \hline
convergence order & &   2.00    &   1.00  &   1.00   \\  \hline
\end{tabular}}
\caption{Acoustic test case $\, \zeta = 10^{-5}$. Errors for density and momentum as the mesh size tends to zero.}
\label{tab-erreurs-acoustique-10-5} \end{table} 

\bigskip \bigskip    \noindent {\bf \large    6) \quad  {Conclusion}}


\smallskip \noindent 
In this work, we studied a lattice Boltzmann scheme based on a geometry of equilateral triangles.
The D2T4 scheme places the physical degrees of freedom at the centers of the triangles.
Since particle trajectories cannot follow straight lines, the classical analysis of lattice Boltzmann
scheme  must be adapted, otherwise incorrect results may emerge.
We conducted this study and the important point is to consider the union of two neighboring
triangles as the basic cell of the model.
We applied this analysis to two model problems in mathematical physics: diffusion and linear acoustics.
To our knowledge, this is the first time that this four-neighbor scheme has been used for acoustics. 

\smallskip \noindent 
The second-order partial differential equations were compared with the results of the schemes
on the one hand and with an analytical calculation on the other. We established the consistency of
the two approaches: the error decreases as the mesh size approaches zero. 

\smallskip \noindent
However, the convergence of the acoustic D2T4 model towards the analytical solution could be improved,
as it is limited to the first order of convergence for both components of the impulsion.
Moreover, the analysis proposed here does not seem entirely satisfactory to us.
In particular, we need to better understand why the velocity distribution 
of the relation~(\ref{w.nabla})
is emerging for the propagation of a bipoint.

\newpage 
\bigskip \bigskip    \noindent {\bf \large   {Acknowledgements}}

This work was initiated in June 2016 and May 2017 when the authors benefited from a stay at the
Beijing Computational Science Research Center. The authors would like to thank the CSRC,
and in particular Professor Li-Shi Luo, for their kind invitation and warm welcome during these two stays. 
FD would like to express his warmest thanks to his partner Sophie Mougel,
who enabled him to complete this work, begun many years ago, thanks to a stay
in La Bresse in the Vosges in February 2026.

\bigskip \bigskip    \noindent {\bf \large    Appendix A}


\smallskip \noindent 
If we duplicate the analysis performed for square D2Q9 for a triangle, for example,
without taking into account the difference between incoming and outgoing particles, Boltzmann's scheme
on a lattice is written as (\ref{evolution-du22}). For a left-type triangle, we have
\moneqstar
\Lambda_\ell = M^\ell \, \, {\rm diag} (v^\ell . \nabla)  \,\, M^{-\ell}
\monendstar
with the left velocities $ \, v^\ell \, $ given by the relation (\ref{vitesses-gauche}) 
and the particles to moments matrix $ \, M^\ell \, $ by the relation (\ref{matrice-M-gauche}).
Then we have 
%
%
\moneq \label{Lambda-naive-gauche}
\Lambda_\ell = \begin{pmatrix} 0 & \partial_x &  \partial_y & 0 \\
{3\over8} \lambda^2 \partial_x & {1\over2} \lambda \, \partial_x & -{1\over2} \lambda \, \partial_y & {1\over8} \, \partial_x  \\ 
{3\over8} \lambda^2 \partial_y & -{1\over2} \lambda \, \partial_y & -{1\over2} \lambda \, \partial_x & {1\over8} \, \partial_y \\
0 & {1\over2} \lambda^2 \, \partial_x & {1\over2} \lambda^2 \, \partial_y & 0 
\end{pmatrix} \, . 
\monend
For a single conservation law, the decomposition
\moneqstar
\Lambda  =  \begin{pmatrix} A & B \\ C & D \end{pmatrix} 
\monendstar
introduces the following matrices
\moneqstar \left\{    \begin{array} {rl} 
A = \big( 0 \big) \,,\,\, & B =  \begin{pmatrix}  \partial_x &  \partial_y & 0 \end{pmatrix} \\
C =   \begin{pmatrix} {3\over8} \lambda^2 \partial_x \\ {3\over8} \lambda^2 \partial_y \\ 0  \end{pmatrix} ,\, &
D =  \begin{pmatrix}  {1\over2} \lambda \, \partial_x & -{1\over2} \lambda \, \partial_y & {1\over8} \, \partial_x  \\
  -{1\over2} \lambda \, \partial_y & -{1\over2} \lambda \, \partial_x & {1\over8} \, \partial_y \\
   {1\over2} \lambda^2 \, \partial_x & {1\over2} \lambda^2 \, \partial_y & 0 \end{pmatrix} \, .  \end{array} \right. 
\monendstar
With the equilibrium (\ref{equilibre-thermique}), we have
\moneqstar
Y^{\rm eq}  =  E \, ( \rho ) \,,\,\, E =  \begin{pmatrix} 0 \\ 0 \\ \alpha \, \lambda^2   \end{pmatrix}  \,. 
\monendstar
If we apply the Taylor expansion method in the ABCD framework  \cite{Du22}
for the discrete evolution (\ref{evolution-du22}).
We get at fourth order accuracy the partial differential equation
\moneq \label{equivalent-edp-scalaire-ordre-4}
\partial_{t} \rho + \alpha_1 \, \rho + \Delta t \, \alpha_2 \, \rho + \Delta t^2 \, \alpha_3 \, \rho
+ \Delta t^3 \, \alpha_4 \, \rho =  {\rm O} (\Delta t^4)
\monend
with operators $ \, \alpha_j \, $ obtained by the ``Berlin'' algorithm \cite{ADGL14}: 
\moneq \label{algorithme-berlin-ordre-4} \left\{ \begin{array}   {l}
\alpha_1 = A + B \, E  \\
\beta_1 =  E \, \alpha_1 - (C +  D \, E)   \\ 
\alpha_2 =  B \,  \Sigma \, \beta_1   \\ 
\beta_2 = \Sigma \, \beta_1 \, \alpha_1 + E \, \alpha_2 - D \, \Sigma \, \beta_1 \\ 
\alpha_3 = B \, \Sigma \, \beta_2 + {1\over12}\, B_2 \, \beta_1 - {1\over6} B \, \beta_1 \, \alpha_1 \\ 
\beta_3 = \Sigma \, \beta_1\, \alpha_2 + E \, \alpha_3 - D\, \Sigma\, \beta_2 + \Sigma \, \beta_2 \, \alpha_1
+ {1\over6} \, D \, \beta_1 \, \alpha_1 - {1\over12} \, D_2 \, \beta_1 - {1\over12} \, \beta_1 \, \alpha_1^2  \\
\alpha_4 = B \, \Sigma \, \beta_3 + {1\over4} \, B_2 \, \beta_1 +  {1\over6} \, B \, D_2 \, \Sigma \, \beta_1 -  {1\over6} \, A \, B \, \beta_2 \\ 
\qquad \quad
-  {1\over6} \, B \, E \, \alpha_1 \, \alpha_2  -  {1\over6} \, B \, E \, \alpha_2 \, \alpha_1   -  {1\over6} \, B \, \Sigma \, \alpha_1^2  \,,
\end{array} \right. \monend
%
%
with 
\moneqstar
\Lambda^2 \equiv  \begin{pmatrix} A_2 & B_2 \\ C_2 & D_2 \end{pmatrix} \, . 
\monendstar
With the help of formal calculus \cite{sage}, we obtain without difficulty
\moneq  \left\{ \begin{array} {l} \label{edp-scalaire-2011-gauche}
\alpha_1^\ell = 0 \\
\alpha_2^\ell =   -{1\over8} \, \lambda^2 \, (\alpha+3) \, \sigma_j \, (\partial_x^2 + \partial_y^2 ) \\
\alpha_3^\ell = {1\over192} \, \lambda^3 \, (12 \, \sigma_j^2-1) \,  (\alpha+3) \, \partial_x \, \big(\partial_x^2 - 3 \, \partial_y^2 \big) \\
\alpha_4^\ell = {1\over256} \,  \lambda^4 \,  \,  (\alpha+3) \, \sigma_j \,
\big( 8 \, \sigma_j^2 + (1-\alpha)\,(1 - 4 \, \sigma_j \, (\sigma_j + \sigma_e) \big) \, (\partial_x^2 + \partial_y^2)^2
\end{array} \right. \monend
%
%
with the anisotropic operator $ \, \partial_x^2 - 3 \partial_y^2 \, $ at third order. 
The result is consistent with the contribution \cite{DL13} (formula (6.2)). 

\bigskip \monitem
For a right-type triangle, we replace the operator matrix $ \, \Lambda_\ell \, $ by $ \, \Lambda_r \, $
now defined by 
\moneqstar
\Lambda_r  = M^r  \, \, {\rm diag} (v^r . \nabla )  \,\, M^{-r}
\monendstar
with the right velocities $ \, v^r \, $ given by the relation (\ref{vitesses-droite}) 
and the particles to moments matrix~$ \, M^r \, $ by the relation (\ref{matrice-M-droite}).
Then we have 
\moneq \label{Lambda-naive-droite}
\Lambda_r = \begin{pmatrix} 0 & \partial_x &  \partial_y & 0 \\
{3\over8} \lambda^2 \partial_x & -{1\over2} \lambda \, \partial_x & {1\over2} \lambda \, \partial_y & {1\over8} \, \partial_x  \\ 
{3\over8} \lambda^2 \partial_y & {1\over2} \lambda \, \partial_y & {1\over2} \lambda \, \partial_x & {1\over8} \, \partial_y \\
0 & {1\over2} \lambda^2 \, \partial_x & {1\over2} \lambda^2 \, \partial_y & 0 
\end{pmatrix} \, . 
\monend
%
%
We note that the two advection matrices  (\ref{Lambda-naive-gauche}) and  (\ref{Lambda-naive-droite})   in the moment basis differ. 
This remark has no impact on the implementation of the algorithm (\ref{algorithme-berlin-ordre-4}).
But the results obtained for partial differential equations are somewhat modified for a right-type triangle.
We obtain 

\vskip -.5 cm 
\moneq  \left\{ \begin{array} {l} \label{edp-scalaire-2011-droite}
\alpha_1^r = 0 \\
\alpha_2^r =   -{1\over8} \, \lambda^2 \, (\alpha+3) \, \sigma_j \, (\partial_x^2 + \partial_y^2 ) \\
\alpha_3^r = {1\over192} \, \lambda^3 \, (12 \, \sigma_j^2-1) \,  (\alpha+3) \, \partial_x \, \big(-\partial_x^2 + 3 \, \partial_y^2 \big) \\
\alpha_4^r = {1\over256} \,  \lambda^4 \,  \,  (\alpha+3) \, \sigma_j \,
\big( 8 \, \sigma_j^2 + (1-\alpha)\,(1 - 4 \, \sigma_j \, (\sigma_j + \sigma_e) \big) \, (\partial_x^2 + \partial_y^2)^2 \, . 
\end{array} \right. \monend
Comparing (\ref{edp-scalaire-2011-gauche}) and  (\ref{edp-scalaire-2011-droite}), we have $ \, \alpha_3^\ell + \alpha_3^r = 0 $.
The two systems of equations differ for the third-order operator!

\bigskip \bigskip    \noindent {\bf \large    Appendix B}

\smallskip \noindent
With the same assumptions as in Appendix A, the equivalent equations for the acoustic system use the same matrices
$ \, \Lambda_\ell \, $ and $ \,  \Lambda_r \, $ defined in   (\ref{Lambda-naive-gauche}) and (\ref{Lambda-naive-droite}) 
respectively.
The ABCD structure now uses an A block of order 3. We have for a left-type triangle 
\moneqstar \left\{    \begin{array} {rl} 
A = \begin{pmatrix}     0 & \partial_x &  \partial_y \\
{3\over8} \lambda^2 \partial_x & {1\over2} \lambda \, \partial_x & -{1\over2} \lambda \, \partial_y  \\
{3\over8} \lambda^2 \partial_y & -{1\over2} \lambda \, \partial_y & -{1\over2} \lambda \, \partial_x & {1\over8} \end{pmatrix}
\,,\,\,& B =  \begin{pmatrix} 0 \\  {1\over8} \, \partial_x \\ {1\over8} \, \partial_y  \end{pmatrix} \\
C =  \begin{pmatrix}  0 & {1\over2} \lambda^2 \, \partial_x & {1\over2} \lambda^2 \, \partial_y \end{pmatrix} \,,\,\, & 
D = \big( 0 \big) \, . 
\end{array} \right. \monendstar
The matrix $ \, E \, $ of equilibria admits the form
\moneq \label{EEE-acoustique} 
e^{\rm eq}  =  Y^{\rm eq}  =  E \, \begin{pmatrix}   \rho \\ J_x \\ J_y   \end{pmatrix}  \,,\,\,
E =  \begin{pmatrix} \alpha \, \lambda^2   & 0 & 0 \end{pmatrix}  \,. 
\monend
The first line of  (\ref{algorithme-berlin-ordre-4}) can be written 

\smallskip \noindent
$ \alpha_1 = A + B \, E  =  A + \begin{pmatrix} 0 \\  {1\over8} \, \partial_x \\ {1\over8} \, \partial_y  \end{pmatrix} \,
\begin{pmatrix} \alpha \, \lambda^2  & 0 & 0 \end{pmatrix}  $

\smallskip \noindent 
$ \,\,\, =  \begin{pmatrix}     0 & \partial_x &  \partial_y \\
{3\over8} \lambda^2 \partial_x & {1\over2} \lambda \, \partial_x & -{1\over2} \lambda \, \partial_y  \\
{3\over8} \lambda^2 \partial_y & -{1\over2} \lambda \, \partial_y & -{1\over2} \lambda \, \partial_x \end{pmatrix}
+  \begin{pmatrix} 0 & 0 & 0  \\  {1\over8} \lambda^2 \, \alpha \, \partial_x & 0 & 0 \\
{1\over8} \lambda^2 \, \alpha \, \partial_y & 0 & 0    \end{pmatrix}
= \begin{pmatrix}     0 & \partial_x &  \partial_y \\
{1\over8} \lambda^2 \, (3 + \alpha) \, \partial_x & {1\over2} \lambda \, \partial_x & -{1\over2} \lambda \, \partial_y  \\
{1\over8} \lambda^2 \, (3 + \alpha) \, \partial_y &  -{1\over2} \lambda \, \partial_y & -{1\over2} \lambda \, \partial_x \end{pmatrix} $.

\smallskip \noindent
Then at first order, the acoustic system for left-type triangles is anisotrop and takes the expression
\moneq  \left\{ \begin{array} {rcl} \label{acoustique-gauche-ordre-1}
\partial_t \rho + {\rm div} J &=& {\rm O}(\Delta t) \\
\partial_t  J_x + {1\over8} \lambda^2 \, (3 + \alpha) \, \partial_x \rho + {1\over2} \lambda \,  \big( \partial_x J_x -  \partial_y J_y \big) 
&=& {\rm O}(\Delta t) \\
\partial_t  J_y + {1\over8} \lambda^2 \, (3 + \alpha) \, \partial_y \rho -  {1\over2} \lambda \,  \big( \partial_y J_x +  \partial_y J_y \big) 
&=& {\rm O}(\Delta t) \, .  \end{array} \right. \monend 

\monitem
For the right type triangle, we start from the matrix $ \,  \Lambda_r \, $ introduced
in   (\ref{Lambda-naive-droite}) and we have
\moneqstar \left\{    \begin{array} {rl} 
A = \begin{pmatrix}     0 & \partial_x &  \partial_y \\
{3\over8} \lambda^2 \partial_x & -{1\over2} \lambda \, \partial_x & {1\over2} \lambda \, \partial_y  \\
{3\over8} \lambda^2 \partial_y & {1\over2} \lambda \, \partial_y & {1\over2} \lambda \, \partial_x & {1\over8} \end{pmatrix}
\,,\,\,& B =  \begin{pmatrix} 0 \\  {1\over8} \, \partial_x \\ {1\over8} \, \partial_y  \end{pmatrix} \\
C =  \begin{pmatrix}  0 & {1\over2} \lambda^2 \, \partial_x & {1\over2} \lambda^2 \, \partial_y \end{pmatrix} \,,\,\, & 
D = \big( 0 \big) \, . 
\end{array} \right. \monendstar
The matrix $ \, E \, $ of equilibria is still given by the relation (\ref{EEE-acoustique})
and we have

\smallskip \noindent
$ \alpha_1 = A + B \, E  =
\begin{pmatrix}     0 & \partial_x &  \partial_y \\
{3\over8} \lambda^2 \partial_x & -{1\over2} \lambda \, \partial_x & {1\over2} \lambda \, \partial_y  \\
{3\over8} \lambda^2 \partial_y & {1\over2} \lambda \, \partial_y & {1\over2} \lambda \, \partial_x \end{pmatrix}
+  \begin{pmatrix} 0 & 0 & 0  \\  {1\over8} \lambda^2 \, \alpha \, \partial_x & 0 & 0 \\
  {1\over8} \lambda^2 \, \alpha \, \partial_y & 0 & 0    \end{pmatrix} $

\smallskip \noindent \quad 
$ = \begin{pmatrix}     0 & \partial_x &  \partial_y \\
{1\over8} \lambda^2 \, (3 + \alpha) \, \partial_x & -{1\over2} \lambda \, \partial_x & {1\over2} \lambda \, \partial_y  \\
{1\over8} \lambda^2 \, (3 + \alpha) \, \partial_y & {1\over2} \lambda \, \partial_y & {1\over2} \lambda \, \partial_x \end{pmatrix} $.

\smallskip \noindent 
Then in consequence, the acoustic system for right-type triangles takes the expression
\moneq  \left\{ \begin{array} {rcl} \label{acoustique-droite-ordre-1}
\partial_t \rho + {\rm div} J &=& {\rm O}(\Delta t) \\
\partial_t  J_x + {1\over8} \lambda^2 \, (3 + \alpha) \, \partial_x \rho - {1\over2} \lambda \,  \big( \partial_x J_x -  \partial_y J_y \big) 
&=& {\rm O}(\Delta t) \\
\partial_t  J_y + {1\over8} \lambda^2 \, (3 + \alpha) \, \partial_y \rho +  {1\over2} \lambda \,  \big( \partial_y J_x +  \partial_y J_y \big) 
&=& {\rm O}(\Delta t) \, .  \end{array} \right. \monend 
We obtain a second system of anisotropic equations for acoustics.
The incorrect gradient-type terms in impulse space change sign between equations  (\ref{acoustique-gauche-ordre-1})
and  (\ref{acoustique-droite-ordre-1}). 

\bigskip \bigskip    \noindent {\bf \large    Appendix C}

\smallskip \noindent 
In this appendix, we present an analysis of the D2T4 scheme for acoustics,
which takes into account network breathing but treats triangles in a decoupled manner. 
We suppose here that the relations  (\ref{g-droite-j}) and (\ref{g-gauche-j})
are valid not only for $ \, 1 \leq j \leq 3 \, $  but also for $ \, j=0 $.
We note that this assumption is incorrect because for $ \, j=0 $, the relations
(\ref{g-gauche-zero}) and (\ref{g-droite-zero})  apply.

\bigskip
For a left-type triangle, we suppose that we have

\moneqstar
g_{\ell j}(x^\ell ,\, t+\Delta t) = f_{r j}^* (x_j^\ell ,\, t) =  f_{rj}^* ( x^\ell - v_j^r \, \Delta t  , \, t)
\monendstar
for $ \, 0 \leq j \leq 3 $.  Then
\moneqstar
g_{\ell}(x  ,\, t+\Delta t) = \exp \big( {\rm diag} ( - v^r . \nabla ) \big) \,\,   f_{r}^* (x ,\, t) 
\monendstar
and
\moneqstar
m(x,\, t+\Delta t) = M^r \, g_{\ell}(x^\ell ,\, t+\Delta t)  = M^r \,  \exp \big( {\rm diag} ( - v^r . \nabla ) \big) \,
M^{-r} \,\,  m^* (x ,\, t)  \, .
\monendstar 
In consequence, we have
\moneqstar
m(x,\, t+\Delta t) = \exp(-\Lambda_r \, \Delta t) \, J_0 \,\, m(x,\, t) \, .
\monendstar
The calculation ends as in Appendix B, and the first-order equations exactly compose the relations (\ref{acoustique-droite-ordre-1}).

\bigskip \noindent 

\bigskip
For a right-type triangle, we suppose now that the discrete time iteration 
\moneqstar
g_{r j}(x^r ,\, t+\Delta t) = f_{\ell j}^* (x_j^r ,\, t) =  f_{r\ell j}^* ( x^r - v_j^\ell \, \Delta t  , \, t)
\monendstar
is correct  for $ \, 0 \leq j \leq 3 $.  Then
\moneqstar
g_{r}(x  ,\, t+\Delta t) = \exp \big( {\rm diag} ( - v^\ell . \nabla ) \big) \,\,   f_{\ell}^* (x ,\, t) 
\monendstar
and
\moneqstar
m(x,\, t+\Delta t) = M^\ell \, g_{r}(x^r ,\, t+\Delta t)  = M^\ell \,  \exp \big( {\rm diag} ( - v^\ell . \nabla ) \big) \,
M^{-\ell} \,\,  m^* (x ,\, t)  \, .
\monendstar
Therefore,
\moneqstar
m(x,\, t+\Delta t) = \exp(-\Lambda_\ell \, \Delta t) \, J_0 \,\, m(x,\, t) \, .
\monendstar
The calculation ends as in Appendix B, and the first-order equations now satisfy the relations~(\ref{acoustique-gauche-ordre-1})!
With one reversal between the left-type and right-type triangles, the conclusions in Appendix B remain unchanged. 
%

\bigskip \bigskip      \noindent {\bf  \large  References }


\end{document}